\documentclass[a4paper,11pt]{amsart}

\newfont{\cyr}{wncyr10 scaled 1100}

\usepackage[left=2.7cm,right=2.7cm,top=3.5cm,bottom=3cm]{geometry}
\usepackage{amsthm,amssymb,amsmath,amsfonts,mathrsfs,amscd}
\usepackage[latin1]{inputenc}
\usepackage[all]{xy}
\usepackage{latexsym}
\usepackage{longtable}

\theoremstyle{plain}
\newtheorem{theorem}{Theorem}[section]
\newtheorem{corollary}[theorem]{Corollary}
\newtheorem{lemma}[theorem]{Lemma}
\newtheorem{proposition}[theorem]{Proposition}

\theoremstyle{definition}
\newtheorem{definition}[theorem]{Definition}

\theoremstyle{remark}
\newtheorem{notation}[theorem]{Notation}
\newtheorem{remark}[theorem]{Remark}

\newcommand{\Emb}{\operatorname{Emb}}

\newcommand{\E}{\mathcal E}
\newcommand{\g}{\gamma}
\newcommand{\G}{\Gamma}

\newcommand{\Q}{\mathbb{Q}}
\newcommand{\Z}{\mathbb{Z}}

\newcommand{\C}{\mathbb{C}}

\newcommand{\PP}{\mathbb{P}}

\newcommand{\Gal}{\operatorname{Gal}}
\newcommand{\GL}{\operatorname{GL}}

\newcommand{\Div}{\operatorname{Div}}

\newcommand{\ord}{{\operatorname{ord}}}
\newfont{\gotip}{eufb10 at 12pt}

\newcommand{\cO}{{\mathcal O}}

\newcommand{\om}{{\omega}}

\newcommand{\SL}{{\mathrm {SL}}}

\newcommand{\Pic}{{\mathrm{Pic}}}

\newcommand{\R}{{\mathbb R}}
\newcommand{\M}{{\mathrm{M}}}

\newcommand{\PGL}{{\mathrm{PGL}}}

\newcommand{\W}{\mathbb W}

\def \mint {\times \hskip -1.1em \int}

\DeclareMathOperator{\Hom}{Hom} 

\newcommand{\X}{\mathbb X}
\newcommand{\Y}{\mathbb Y}
\newcommand{\U}{\mathbb U}

\newcommand{\res}{\mathrm{res}}
\newcommand{\D}{\mathbb D}

\newcommand{\fr}{\mathfrak}
\newcommand{\cl }{\mathcal}

\newcommand{\longmono}{\mbox{$\lhook\joinrel\longrightarrow$}}

\newcommand{\smallmat}[4]{\bigl(\begin{smallmatrix}#1&#2\\#3&#4\end{smallmatrix}\bigr)}

\include{thebibliography}

\begin{document}

\title[Quaternionic Darmon points over genus fields]{The rationality of quaternionic Darmon points over genus fields of real quadratic fields}
\author{Matteo Longo and Stefano Vigni}

\begin{abstract} 
Darmon points on $p$-adic tori and Jacobians of Shimura curves over $\Q$ were introduced in \cite{LRV1} and \cite{LRV2} as generalizations of Darmon's Stark--Heegner points. In this article we study the algebraicity over extensions of a real quadratic field $K$ of the projections of Darmon points to elliptic curves. More precisely, we prove that linear combinations of Darmon points on elliptic curves weighted by certain genus characters of $K$ are rational over the predicted genus fields of $K$. This extends to an arbitrary quaternionic setting the main theorem on the rationality of Stark--Heegner points obtained by Bertolini and Darmon in \cite{BD-annals}, and at the same time gives evidence for the rationality conjectures formulated in the joint paper with Rotger \cite{LRV2} and by M. Greenberg in \cite{Gr}. In light of this result, quaternionic Darmon points represent the first instance of a systematic supply of points of Stark--Heegner type other than Darmon's original ones for which explicit rationality results are known. 
\end{abstract}

\address{Dipartimento di Matematica Pura e Applicata, Universit\`a di Padova, Via Trieste 63, 35121 Padova, Italy}
\email{mlongo@math.unipd.it}
\address{Departament de Matem\`atica Aplicada II, Universitat Polit\`ecnica de Catalunya, C. Jordi Girona 1-3, 08034 Barcelona, Spain}
\email{stefano.vigni@upc.edu}



\maketitle

\section{Introduction} \label{intro}

Given an elliptic curve $E$ over $\Q$ of conductor $N$, a prime number $p$ dividing $N$ exactly and a real quadratic field $K$ in which $p$ is inert and all primes dividing $N/p$ split (Stark--Heegner hypothesis), Darmon introduced in \cite{Dar} the notion of Stark--Heegner points on $E$. These local points live in $E(K_p)$, where $K_p$ is the completion of $K$ at $p$, but are conjectured to be rational over (narrow) ring class fields of $K$ and to satisfy a Shimura reciprocity law under Galois actions. In fact, the arithmetic properties of Stark--Heegner points are expected to be as rich as those enjoyed by classical Heegner points, which are defined via the theory of complex multiplication and are known to be rational over ring class fields of imaginary quadratic fields. 

Stark--Heegner points were later lifted to certain quotients of classical modular Jacobians by Dasgupta in \cite{Das}; this was done by proving a rigid analytic uniformization result for modular Jacobians which boils down to an equality of $\mathcal L$-invariants and turns out to be a strong form of a theorem of Greenberg and Stevens (\cite{GS}). It is important to observe that both Darmon's and Dasgupta's strategies rely heavily on the theory of modular symbols (and so on the presence of cusps on classical modular curves), and thus they do not extend to more general situations where the Stark--Heegner condition is not satisfied but considerations on the signs of the functional equations of the relevant $L$-functions predict the existence of points of Stark--Heegner type. 

To overcome this difficulty, in \cite{Gr} M. Greenberg replaced modular symbols by a general cohomological study of Shimura varieties over totally real fields and reinterpreted Darmon's theory in terms of group cohomology. The outcome of his work is, in the simplest case where the base field is $\Q$, a conjectural construction of local Stark--Heegner points on $E$ which are expected to be rational over appropriate ring class fields of a real quadratic field $K$. In this context, the field $K$ satisfies a modified Stark--Heegner hypothesis forcing the number of primes dividing $N/p$ and inert in $K$ to be even. 

In \cite{LRV1} we gave an explicit rigid analytic uniformization of the maximal toric quotient of the Jacobian of a Shimura curve over $\Q$ at a prime dividing exactly the level, a result that can be viewed as complementary to the classical theorem of \v{C}erednik and Drinfeld which provides rigid analytic uniformizations at primes dividing the discriminant (\cite{BC}). More precisely, denote $J_0^D(Mp)$ the Jacobian variety of the Shimura curve of discriminant $D>1$ and level $Mp$. Combining techniques in group cohomology and $p$-adic integration in the spirit of the works of Dasgupta and Greenberg, we exhibited a $p$-adic torus which is isogenous to the product of two copies of the $p$-new quotient of $J_0^D(Mp)$. As a corollary, we offered a proof of the isogeny conjecture formulated in \cite[Conjecture 2]{Gr} (an alternative proof of which has been given by Dasgupta and Greenberg in \cite{DG}): this made Greenberg's construction of local points on elliptic curves over $\Q$ unconditional, but left the rationality conjecture \cite[Conjecture 3]{Gr} wide open. 
 
As an application of our work on rigid analytic uniformizations, in \cite{LRV2} we introduced a new family of Stark--Heegner type points on $p$-adic tori and Jacobians of Shimura curves, which we called \emph{Darmon points}. These are, at the same time, lifts of Greenberg's points on elliptic curves and generalizations of Dasgupta's Stark--Heegner points on modular Jacobians, and we formulated conjectures about their rationality over (narrow) ring class fields of real quadratic fields and their Galois properties. Moreover, assuming these rationality conjectures and applying the methods of \cite{LV1}, we proved the conjecture of Birch and Swinnerton-Dyer for elliptic curves over ring class fields of real quadratic fields in the case of analytic rank $0$.  

At present, the central open problems in the theories of Stark--Heegner and Darmon points are the rationality conjectures. The first result towards the rationality of Stark--Heegner points was obtained by Bertolini and Darmon in \cite{BD-annals}, where they proved that linear combinations of Stark--Heegner points weighted by genus characters of $K$ are rational over the predicted genus fields of $K$. To achieve this, the authors showed that linear combinations of this kind can be expressed in terms of Heegner points and applied their previous results on Hida families of modular forms (\cite{BD}).

In this paper we prove the analogue of the results of \cite{BD-annals} for projections to elliptic curves of the quaternionic Darmon points introduced in \cite{LRV2}. In other words, we show that linear combinations with genus character coefficients of Darmon points on elliptic curves are rational over the fields that were predicted in \cite{LRV2}. In light of this feature, quaternionic Darmon points represent the first instance of a systematic supply of points of Stark--Heegner type other than Darmon's original ones for which explicit rationality results are known. Thus our main theorem provides some evidence for both the rationality conjecture formulated in \cite{LRV2} and the conjecture proposed by Greenberg in \cite[Conjecture 3]{Gr} when the totally real base field is $\Q$. The present article should be viewed as a natural step further in the line of investigation begun in \cite{LRV1} and carried on in \cite{LRV2}.  

While the reader familiar with the work of Bertolini and Darmon will notice that our strategy is inspired by \cite{BD-annals} (e.g., the main result of \cite{BD} plays a crucial role in the last part of our arguments), the methods we use are of a different nature and build on the constructions of \cite{LRV1} in an essential way. In fact, in this article we develop a series of new techniques and obtain various auxiliary results which may be of independent interest. Finally, we remark that, as in \cite{BD-annals}, the scheme of proof presents formal analogies with that of Kronecker's solution to Pell's equation (see \cite[Introduction]{BD-annals} for an illustration of this parallelism). 

Now we describe our results more in detail. Fix an elliptic curve $E$ over $\Q$ of square-free conductor $N$ and a real quadratic field $K$ of discriminant $d_K$ prime to $N$, whose ring of integers we denote $\cO_K$. Suppose that there is a factorization $N=MDp$ such that 
\begin{itemize}
\item all the primes dividing $M$ (respectively, $Dp$) are split (respectively, inert) in $K$; 
\item the number of primes dividing $D$ is \emph{even} and $\geq2$.
\end{itemize} 
Let $B$ denote the indefinite quaternion algebra over $\Q$ of discriminant $D$ and choose Eichler orders $R_0^D(Mp)\subset R_0^D(M)$ of $B$ of level $Mp$ and $M$, respectively. Write $\Gamma_0^D(Mp)$ and $\Gamma_0^D(M)$ for the groups of elements of norm $1$ in $R_0^D(Mp)$ and $R_0^D(M)$, respectively, and consider the Ihara group 
\[ \Gamma:=\Big\{\gamma\in R_0^D(M)\otimes_\Z\Z[1/p]\;\Big|\;{\rm norm}(\gamma)=1\Big\}. \]
For any sign $\pm$ fix an isomorphism $H^1(E(\C),\Z)^\pm\simeq\Z$, where the superscript denotes the $\pm$-eigenspace for complex conjugation. Following \cite[\S 7.7]{LRV1}, we define a canonical element 
\[ \boldsymbol{\mu}_E^\pm\in H^1\bigl(\Gamma,\mathcal M_0(\Q)\bigr) \] 
where $\mathcal M_0(\Q)$ is the $\Q$-vector space of $\Q$-valued measures on $\PP^1(\Q_p)$ with total mass $0$. Its image in $H^1\bigl(\Gamma_0^D(Mp),\Q\bigr)$, obtained by applying Shapiro's lemma, spans the one-dimensional $\Q$-vector space on which the archimedean involution $W_\infty$ acts via $\pm$ and the Hecke algebra of level $N$ acts via the character associated with the newform $f$ attached to $E$ by modularity (see \S \ref{Hecke+} for details on these actions). 
  
Fix a sign $\pm$. Using techniques introduced in \cite{LRV1} and \cite{LRV2}, in \S \ref{sec2.5} we use $\boldsymbol\mu_E^\pm$ to attach a Darmon point $P_\psi^\pm\in E(K_p)$ to (the $\Gamma_0^D(M)$-conjugacy class of) each optimal embedding $\psi$ of $\cO_K$ into $R_0^D(M)$: these points are the object of our investigation in the present paper. 

Let $G_K^+$ be the narrow class group of $K$ (i.e., the Galois group over $K$ of the narrow Hilbert class field of $K$). For a genus character $\chi$ of $K$  (i.e., an unramified quadratic character of $G_K^+$) we define the point  
\[ P_\chi:=\sum_{\sigma\in G_K^+}P_{\psi_\sigma}^\epsilon\in E(K_p), \] 
where the Galois action $\psi\mapsto\psi_\sigma$ on (conjugacy classes of) optimal embeddings is described in \S \ref{Galois-emb-subsec} and the choice of $\epsilon\in\{\pm\}$ depends on $\chi$. 

The main result of this paper, which is Theorem \ref{theorem5.1} in the text, is the following 
 
\begin{theorem} \label{main-thm-intro}
Let $\chi$ be a genus character of $K$ associated with a pair $(\chi_1,\chi_2)$ of Dirichlet characters such that $\chi_i(-MD)=-w_{MD}$ for $i=1,2$, where $w_{MD}$ is the eigenvalue of the Atkin--Lehner involution $W_{MD}$ acting on $f$.   
\begin{enumerate}
\item There exists $n\in\Z$ such that $nP_\chi\in E(H_\chi)$ where $H_\chi$ is the genus field of $K$ cut out by $\chi$. 
\item The point $nP_\chi$ is of infinite order if and only if $L'(E/K,\chi,1)\not=0$.
\end{enumerate} 
\end{theorem} 

Now we outline the strategy of proof of this theorem. In doing this, we state two results (Theorems \ref{intro-thm-1} and \ref{thm-intro-3}) which may be of independent interest. 

Write $\D$ for the module of $\Z_p$-valued measures on $\Y:=\Z_p^2$ which are supported on the subset $\X$ of primitive elements (i.e., those $(x,y)\in\Y$ which are not divisible by $p$). In Section \ref{sec-control} we exploit a control theorem for Hida families in the quaternionic setting (analogous to the one proved by Greenberg and Stevens in \cite{GS}) to construct a canonical element 
\[ \boldsymbol{\tilde\mu}_E^\pm\in H^1\bigl(\Gamma_0^D(M),\D\bigr)\otimes_{\Z_p}\Q_p \] 
such that the image of $\boldsymbol{\tilde\mu}_E^\pm$ via the map $\pi:\X\rightarrow\PP^1(\Q_p)$ taking $(x,y)$ to $x/y$ coincides with $\boldsymbol\mu_E^\pm$. We fix a representative $\tilde\mu_{E}^\pm$ of this class, which can be chosen so that $\pi_*(\tilde\mu_E^\pm)=\mu_E^\pm$. In order to describe $\boldsymbol{\tilde\mu}_E^\pm$ more precisely, embed $\Z$ into the weight space
\begin{equation} \label{X}
\mathcal X:=\Hom(\Z_p^\times,\Z_p^\times)\simeq\Z/(p-1)\Z\times(1+p\Z_p) 
\end{equation}
by sending $k$ to the map $\bigl(x\mapsto x^{k-2}\bigr)$. Hida theory associates with $f$ a neighbourhood $\U$ of $2$ in $\mathcal X$ (which we assume small enough so that all the results we will state below hold for $\U$) and a formal $q$-expansion 
\begin{equation} \label{q-exp}
f_\infty=\sum_{n=1}^\infty a_n(\kappa)q^n
\end{equation} 
where $a_1=1$ and $a_n$ is a rigid analytic function on $\U$ such that $\sum_{n=1}^\infty a_n(k)q^n$ for an even integer $k\geq2$ is the $q$-expansion of a $p$-stabilized weight $k$ eigenform $f_k$ on $\Gamma_0(N)$ and $f_2=f$. The form $f_k$ is associated, via the Jacquet--Langlands correspondence, with  a modular form $f_k^{\rm JL}$ of weight $k$ and level structure $\Gamma_0^D(Mp)$ (so $f_k^{\rm JL}$ is well defined only up to scalars). Thanks to a result of Matsushima and Shimura (\cite{MS}), $f_k^{\rm JL}$ corresponds to an element $\boldsymbol\phi_k$ in the first cohomology group of $\Gamma_0^D(Mp)$ with values in the dual of the $\C$-vector space of homogeneous polynomials of degree $k-2$. Then the property characterizing $\boldsymbol{\tilde\mu}_E^\pm$ asserts that integrating homogeneous polynomials of degree $k-2$ against $\tilde\mu_{E,\gamma}^\pm$ over $\Z_p\times\Z_p^\times$ (with $\gamma\in\Gamma_0^D(Mp)$) defines a multiple of the projection of $\boldsymbol\phi_k$ to the $\pm$-eigenspace for $W_\infty$. See \S \ref{subsec-control} for details.  

Now let 
\begin{equation} \label{Tate}
\Phi_{\rm Tate}:\bar K_p/\langle q\rangle\longrightarrow E(\bar K_p)
\end{equation} 
denote Tate's $p$-adic uniformization, where $q\in p\Z_p$ is Tate's period for $E$ at $p$ and $\bar K_p$ is an algebraic closure of $K_p$. Let $\log_q$ be the branch of the $p$-adic logarithm satisfying $\log_q(q)=0$ and define $\log_E(P):=\log_q\bigl(\Phi_{\rm Tate}^{-1}(P)\bigr)$ for all $P\in E(K_p)$. 

As a piece of notation, let $z_\psi$ denote the fixed point of $\psi(K^\times)$ acting on $\PP^1(\bar\Q_p)$ such that $\psi(\alpha)(z_\psi,1)=\alpha(z_\psi,1)$ for all $\alpha\in K^\times$ (see \S \ref{sec4.1} for precise definitions). Let also $\varepsilon_K$ be a generator of the group of norm $1$ elements in $\cO_K^\times$ such that $\varepsilon_K>1$ with respect to a fixed embedding $K\hookrightarrow\R$ and define $\gamma_\psi:=\psi(\varepsilon_K)$. 

The first auxiliary result proved in the paper is the following 
\begin{theorem} \label{intro-thm-1} 
There is an integer $m\not=0$ such that 
\[ \log_E(P_\psi^\pm)=m\cdot\int_\X\log_q(x-z_\psi y)d{{\tilde\mu}}_{{E,\gamma_\psi}}^\pm. \]
\end{theorem} 

This theorem, which corresponds to Corollary \ref{coro}, allows us to explicitly compute (at least in principle) Darmon points on elliptic curves. Actually, in Theorem \ref{thm2.11} we prove a more general version for Darmon points on the $p$-adic torus $T/L$ by means of which we obtained in \cite{LRV1} a parametrization of $E$ (cf. \S \ref{sec2.2}). 

\begin{remark}
The integer $m$ appearing in Theorem \ref{intro-thm-1} can be made explicit in terms of the exponent of the abelianization of $\Gamma$ and the degree of an isogeny between the torus $T/L$ and the elliptic curve $E$ (see \S \ref{sec2.2}). 
\end{remark} 

In Section \ref{L-function} we use the cohomology class $\boldsymbol{\tilde\mu}_E^\pm$ to introduce a $p$-adic $L$-function $L_p(f_\infty/K,\chi,k)$ attached to $f_\infty$ and  a genus character $\chi$ of $K$, whose $p$-adic variable $k$ ranges in $\U$. The main result proved by Popa in \cite{Po} can then be combined with the interpolation properties of $\boldsymbol{\tilde\mu}_E^\pm$ to show that, for a classical modular form $f_k$ of even weight $k\geq4$ and trivial character in the Hida family, the value of $L_p(f_\infty/K,\chi,k)$ is equal, up to multiplication by an explicit non-zero constant, to the special value of the complex $L$-function of $f_k^\sharp$ twisted by $\chi$, where the form $f_k^\sharp$ is the $p$-stabilization of $f_k$ (see Theorem \ref{teo4.17}).  

Now recall that a genus character $\chi$ of $K$ corresponds to an unordered pair of quadratic Dirichlet characters $\chi_1$ and $\chi_2$ of discriminants $d_1$ and $d_2$, respectively, such that $d_K=d_1d_2$, and for $j=1,2$ let $L_p(f_\infty,\chi_j,k,s)$ denote the Mazur--Kitagawa two-variable $p$-adic $L$-function associated with $f_\infty$ and $\chi_j$. 

Our last result we highlight in this introduction is the following factorization formula, corresponding to Theorem \ref{thm-factor} in the main body of the article.  

\begin{theorem}\label{thm-intro-3} Let $\chi$ be a genus character of $K$ and let $(\chi_1,\chi_2)$ be the associated pair of Dirichlet characters. There exists a $p$-adic analityc function $\eta$ on $\U$ such that $\eta(k)\neq0$ and   
\[ L_p(f_\infty/K,\chi,k)=\eta(k)L_p(f_\infty,\chi_1,k,k/2)L_p(f_\infty,\chi_2,k,k/2) \] 
for all $k\in \U$. Moreover, $\eta(2)\in(\Q^\times)^2$.   
\end{theorem}

Here is a brief explaination of how Theorems \ref{intro-thm-1} and \ref{thm-intro-3} are used to prove Theorem \ref{main-thm-intro}.  
\begin{enumerate}
\item Theorem \ref{thm-intro-3} provides a link between the values at $k=2$ of the second derivatives of $L_p(f_\infty/K,\chi,k)$ and of the restriction of $L_p(f_\infty,\chi_1,k,s)$ to the line $s=k/2$ (here we need to order the pair $(\chi_1,\chi_2)$ so that the sign of the functional equation of the $L$-function $L(E,\chi_1,s)$ is $-1$, which can be done thanks to our assumptions). 
\item On the one hand, the second derivative of $L_p(f_\infty/K,\chi,k)$ evaluated at $k=2$ is essentially computed by the integrals appearing on the right hand side of Theorem \ref{intro-thm-1}, thus providing a link with Darmon points (this is proved in Theorem \ref{thm-derivatives}). 
\item On the other hand, \cite[Theorem 5.4]{BD} shows that the value at $k=2$ of the second derivative of $L_p(f_\infty,\chi_1,k,k/2)$ encodes the logarithm of certain linear combinations of classical Heegner points.
\end{enumerate}  
The three observations above establish a relation between Heegner points and Darmon points, yielding the desired rationality result for $P_\chi$. 
 
\begin{remark} 
The hypothesis that $M$ be square-free is imposed only to be able to apply the results of \cite{Po}, which are proved under this assumption. Granting an extension of \cite{Po} to more general situations, our arguments work equally well under the weaker condition that $M$ be a product of odd powers of distinct primes (and this restriction comes into play exclusively in \S \ref{atkin-genus-subsec} for some calculations with local Hilbert symbols).
\end{remark} 

\section{Darmon points on $p$-adic tori and elliptic curves} 

\subsection{Quaternion algebras and Shimura curves} \label{sec2.1}

Fix throughout the paper embeddings 
\begin{equation} \label{embeddings}
\bar\Q\;\longmono\;\bar\Q_p,\qquad \bar\Q_p\;\longmono\;\C.
\end{equation}
Let $D>1$ be a square-free product of an \emph{even} number of primes and let $M\geq1$ be a square-free integer prime to $D$. We also fix a prime number $p$ not dividing $6MD$. Write $B$ for the (indefinite) quaternion algebra over $\Q$ of discriminant $D$. 

Let $K=\Q(\sqrt{d_K})$ be a real quadratic field with fundamental discriminant $d_K$ such that all the primes dividing $Dp$ are inert in $K$ and all the primes dividing $M$ are split in $K$, and fix an embedding $K\hookrightarrow\R$. For the rest of the paper, fix an auxiliary real quadratic field $F$ such that all the primes dividing $D$ are inert in $F$ while $p$ is split in $F$. In particular, $F\not=K$. The condition on the primes dividing $D$ implies that $F$ is a splitting field for $B$, so we can (and do) choose an isomorphism $i_F:B\otimes_\Q F\overset\simeq\longrightarrow\M_2(F)$ of $F$-algebras. Since $p$ splits in $F$, the embedding $\bar\Q\hookrightarrow\bar\Q_p$ fixed in \eqref{embeddings} and the isomorphism $i_F$ induce an isomorphism $i_p:B\otimes_\Q\Q_p\overset\simeq\longrightarrow\M_2(\Q_p)$ of $\Q_p$-algebras. Moreover, fix an embedding $F\hookrightarrow\R$; this induces via $i_F$ an isomorphism $i_\infty:B\otimes_\Q\R\overset\simeq\longrightarrow\M_2(\R)$ of $\R$-algebras. Finally, for all primes $\ell\nmid Dp$ fix isomorphisms $i_\ell:B\otimes_\Q\Q_\ell\overset\simeq\longrightarrow\M_2(\Q_\ell)$ of $\Q_\ell$-algebras.

\begin{remark}
When describing local arguments at $p$, the fact that $p$ splits in $F$ allows us to work with $\Q_p$ and not with quadratic extensions of it, which simplifies the computations.
\end{remark} 

Let $R^{\rm max}$ be the maximal order of $B$ satisfying the following condition: if $\ell$ is a prime such that $\ell\nmid D$ then $i_\ell(R^{\rm max}\otimes\Z_\ell)=\M_2(\Z_\ell)$. We fix an Eichler order $R_0^D(M)$ of level $M$ contained in $R^{\rm max}$ by requiring that for all primes $\ell|M$ the order $i_\ell\big(R_0^D(M)\otimes\Z_\ell\big)$ is equal to the order $R_0^{\rm loc}(\ell)$ of $\M_2(\Z_\ell)$ consisting of those matrices $\smallmat abcd$ with $c\equiv 0\pmod{\ell}$. Furthermore, we also fix an Eichler order $R_0^D(Mp)\subset R_0^D(M)$ of level $Mp$ by requiring that, in addition, $i_p\bigl(R_0^D(Mp)\otimes\Z_p\bigr)$ is equal to the order $R_0^{\rm loc}(p)$ of $\M_2(\Z_p)$ consisting of those matrices $\smallmat abcd$ with $c\equiv 0\pmod p$. We denote by $\Gamma_0^D(M)$ and $\Gamma_0^D(Mp)$ the groups of elements in $R_0^D(M)$ and $R_0^D(Mp)$, respectively, of reduced norm $1$. 

Finally, let $X_0^D(M)$ (respectively, $X_0^D(Mp)$) denote the (compact) Shimura curve whose corresponding Riemann surface is defined as the analytic quotient $\Gamma_0^D(M)\backslash \mathcal H$ (respectively, $\Gamma_0^D(Mp)\backslash \mathcal H$), where $\mathcal H$ is the complex upper half plane and the elements of $B^\times$ with positive norm act on it by M\"obius (i.e., fractional linear) transformations via $i_\infty$. 

\subsection{Hecke algebras} \label{Hecke+}

For any subgroup $G$ of $B^\times$ consisting of elements having reduced norm $1$ and any subsemigroup $S$ of $B^\times$ such that $(G,S)$ is a Hecke pair in the sense of \cite[\S 1.1]{AS} we denote by $\mathcal H(G,S)$ the Hecke algebra (over $\Z$) of the pair $(G,S)$ whose elements are combinations with integer coefficients 
of double cosets $T(s):=GsG$ for $s\in S$. 

Let $\rm norm$ be the reduced norm of $B$, let $g\mapsto g^*:={\rm norm}(g)g^{-1}$ be the main involution of $B^\times$ and for any subset $S$ of $B^\times$ let $S^*$ denote the image of $S$ under the main involution.
 
If $M$ is a left $\Z[S^*]$-module, the group $H^r(G,M)$ has a natural right action of $\mathcal H(G,S)$ defined on the level of cochains $c\in Z^r(G,M)$ as follows: for $s\in S$ set 
\[ \bigl(c|T(s)\bigr)(\gamma_1,\dots,\gamma_r):=\sum_{i}s_i^*\cdot c\bigl(t_i(\gamma_1),\dots,t_i(\gamma_r)\bigr) \] 
where $GsG=\coprod_{i} Gs_i$ (finite disjoint decomposition) 
and $t_i:G\rightarrow G$ is defined by the equations $Gs_i\gamma=Gs_j$ (for some $j$) and $g_i\gamma=t_i(\gamma)g_j$. We also define a right action of $\mathcal H(G,S)$ on $H_r(G,M)$ by the formula 
\[ c|T(s):=\sum_{i,j}s_i(m_j)\otimes\bigl[t_i(\gamma_{1,j})|\dots|t_i(\gamma_{r,j})\bigr] \]
for $c=\sum_j{m_j}\otimes[\gamma_{1,j}|\dots|\gamma_{r,j}]\in Z_r(G,M)$.

\begin{remark} \label{rem2.1}
The above formalism of Hecke operators slightly differs from the one adopted in \cite[\S 2.1]{LRV1}, where we considered left actions on homology and cohomology groups instead of right actions.
\end{remark}

For primes $\ell$ let $\Sigma_\ell^{\rm loc}$ denote the semigroup of elements in $R^{\rm max}\otimes\Z_\ell$ with non-zero norm. For every prime $\ell|Mp$ let $\Sigma_0^{\rm loc}(\ell)\subset\Sigma_\ell^{\rm loc}$ be the inverse image under $i_\ell$ of the semigroup of matrices $\smallmat abcd\in\GL_2(\Q_\ell)\cap\M_2(\Z_\ell)$ with $a\in\Z_\ell^\times$ and $c\equiv0\pmod{\ell}$. Then define the semigroup 
\[ S_0^D(Mp):=B^+\cap\Bigg(\prod_{\ell|Mp}\Sigma_0^{\rm loc}(\ell)\times\prod_{\ell\nmid Mp}\Sigma_\ell^{\rm loc}\Bigg) \]
where $B^+$ is the subgroup of elements in $B^\times$ of positive norm. Then we can form the commutative Hecke algebra $\mathcal H(Mp)=\mathcal H\bigl(\Gamma_0^D(Mp),S_0^D(Mp)\bigr)$ generated over $\Z$ by the standard elements $T_n=\sum T(\alpha)$ for $n\geq1$ (the sum being over all the elements $\alpha\in S$ of norm $n$) and $T_{n,n}=T(n)$ for $n\geq1$, $n\nmid MDp$. Analogously, we can define the Hecke algebra $\mathcal H(M)$.

As in \cite{LRV1}, consider the Ihara group 
\[ \Gamma:=\Big\{\gamma\in R_0^D(M)\otimes_\Z\Z[1/p]\;\Big|\;{\rm norm}(\gamma)=1\Big\}\;\longmono\;\SL_2(\Q_p) \] 
where ${\rm norm}:B\rightarrow\Q$ denotes the reduced norm map and the injection is induced by the composition of the canonical injection $B\hookrightarrow B\otimes_\Q\Q_p$ with the isomorphism $i_p$ fixed in \eqref{embeddings}. If 
\[ S_p:=B^+\cap\Bigg(\M_2(\Q_p)\times\prod_{\ell|M}\Sigma_0^{\rm loc}(\ell)\times\prod_{\ell\nmid Mp}\Sigma_\ell^{\rm loc}\Bigg) \]
then we may also consider the Hecke algebra $\mathcal H(\Gamma,S_p)$. See \cite[\S 2.2 and \S 2.3]{LRV1} for details.

Now we introduce Atkin--Lehner involutions. For every prime $\ell|MDp$ choose an element $\omega_\ell\in R_0^D(Mp)$ of norm $\ell$; if $\ell|Mp$ then we can (and do) choose $\omega_\ell$ in such a way that $i_\ell(\omega_\ell)=\smallmat 0{-1}{\ell}0$ up to an element of $R_0^{\rm loc}(\ell)$ of norm $1$. It is known that $\omega_\ell$ normalizes both $\Gamma_0^D(Mp)$ and $\Gamma_0^D(M)$ when $\ell|MD$ (cf. \cite[\S 2]{Ogg}), and it turns out that the same is true of $\omega_p$. Finally, choose $\omega_\infty\in R_0^D(Mp)$ of norm $-1$; of course, $\omega_\infty$ normalizes both $\Gamma_0^D(Mp)$ and $\Gamma_0^D(M)$. 

Suppose that a semigroup $S\subset B^\times$ contains $\omega_\infty$ and $\omega_\ell$ for all primes $\ell|MDp$. Then define $W_\infty:=T(\omega_\infty)$ and $W_\ell:=T(\omega_\ell)$ in $\mathcal H(G,S)$; further, for every integer $m|MDp$ set $W_m:=\prod_{\ell|m}W_\ell$. Note that $W_m=T(\omega_m)$ where $\omega_m:=\prod_{\ell|m}\omega_\ell$. In particular, $W_p$ belongs to $\mathcal H(\Gamma,S_p)$ and, in fact, $W_p=U_p$ in this Hecke algebra (see \cite[\S 2.3]{LRV1} for details).

\subsection{Quadratic forms} \label{sec4.1}

Let $\tau$ denote the generator of $\Gal\bigl(F(\sqrt{d_K})/F\bigr)$. If $\psi:K\hookrightarrow B$ is an embedding of $\Q$-algebras then we may consider the quadratic form $Q_\psi(x,y)\in P_2(F)$ associated with $\psi$, which is defined by 
\[ Q_\psi(x,y):=cx^2-2axy-by^2\qquad\text{where $i_F\bigl(\psi(\sqrt{d_K})\bigr)=\begin{pmatrix}a&b\\c&{-a}\end{pmatrix}$}. \] 
We explicitly observe that $c\not=0$. If fact, if $c=0$ then $d=a^2$, whence $K\subset F$, a possibility which is ruled out by our choice of $F$.

We can factor the quadratic form $Q_\psi(x,y)$ as 
\[ Q_\psi(x,y)=c(x-z_\psi y)(x-\bar z_\psi y) \] 
where $z_\psi,\bar z_\psi\in F(\sqrt{d_K})-F$ are the roots of the equation $cz^2-2az-b=0$ and $\tau(z_\psi)=\bar z_\psi$. The two roots $z_\psi$ and $\bar z_\psi$ are the only fixed points for the action of $\psi(K^\times)$ on $\PP^1\big(F(\sqrt{d_K})\big)$ by M\"obius transformations via $i_F$. We may also order $z_\psi$ and $\bar z_\psi$ by requiring that  
\[ \psi(\alpha)\binom{z_\psi}{1}=\alpha\binom{z_\psi}{1},\qquad\psi(\alpha)\binom{\bar z_\psi}{1}=\tau(\alpha)\binom{\bar z_\psi}{1} \]
for all $\alpha\in K$.

Let $\mathcal H_p:=\C_p-\Q_p$ denote Drinfeld's $p$-adic half plane. The completion $K_p$ of $K$ at the prime $(p)$ is the (unique, up to isomorphism) unramified quadratic extension of $\Q_p$. Since $p$ is split in $F$, the completion of $F(\sqrt{d_K})$ with respect to the valuation induced by the embedding $\bar\Q\hookrightarrow\bar\Q_p$ chosen in \eqref{embeddings} is again $K_p$. Thus $z_\psi$ and $\bar z_\psi$ can also be seen as points in $K_p-\Q_p$. Since there are canonical isomorphisms  
\[ \Gal\bigl(F(\sqrt{d_K})/F\bigr)\simeq\Gal(K/\Q)\simeq\Gal(K_p/\Q_p), \] 
in the rest of the article we shall view $\tau$ as the generator of $\Gal(K_p/\Q_p)$ as well.
 
\subsection{Darmon points on tori} \label{sec2.5}\label{sec-2.1} \label{sec-2.2}
 
For any abelian group $A$ denote by $\mathcal M(A)$ the group of $A$-valued measures on $\PP^1(\Q_p)$ and by $\mathcal M_0(A)$ those measures in $\mathcal M(A)$ with total mass $0$. There is a canonical left action of $\GL_2(\Q_p)$ on $\mathcal M(A)$ and $\mathcal M_0(A)$ defined by the integration formula 
\[ \int_{\PP^1(\Q_p)}\varphi(t)d(\gamma\cdot\nu)(t):=\int_{\PP^1(\Q_p)}\varphi\left(\frac{at+b}{ct+d}\right)d\nu \]  
for all step functions $\varphi:\PP^1(\Q_p)\rightarrow\C_p$ and all $\gamma=\smallmat abcd\in\GL_2(\Q_p)$. Then $B^\times$ acts on $\mathcal M(A)$ and $\mathcal M_0(A)$ via $i_p$ and the embedding $B\hookrightarrow B\otimes\Q_p$. 

The group $H^1\bigl(\Gamma,\mathcal M_0(\Q)\bigr)$ is endowed with actions of the Hecke algebra $\mathcal H(\Gamma,S_p)$ considered in \S \ref{Hecke+} and of the involutions $W_m$ and $W_\infty$. For any choice of sign $\pm$, define $H^1\bigl(\Gamma,\mathcal M_0(\Q)\bigr)^{f,\pm}$ to be the subspace of $H^1\bigl(\Gamma,\mathcal M_0(\Q)\bigr)$ consisting of those $\xi$ such that $\xi|T=\theta_f(T)\xi$ for all $T\in\mathcal H(Mp)$ and $\xi|W_\infty=\pm\xi$ where $\theta_f:\mathcal H(Mp)\rightarrow\Z$ is the morphism associated with $f$. By \cite[Proposition 25]{Gr}, the $\Q$-vector space $H^1\bigl(\Gamma,\mathcal M_0(\Q)\bigr)^{f,\pm}$ has dimension $1$. In \cite[Section 4]{LRV1} we constructed an explicit generator for this vector space. Since we will need this description in the following computations, we will briefly review setting and results of \cite{LRV1}. 

Recall the identification between $H^1\bigl(\Gamma,\mathcal M_0(A)\bigr)$ with $H^1\bigl(\Gamma,\mathcal F_{\rm har}(A)\bigr)$ where, for any abelian group $A$, the symbol $\mathcal F_{\rm har}(A)$ stands for the abelian group of $A$-valued harmonic cocycles (see, e.g., \cite[Lemma 27]{Gr} and \cite[\S 4.1]{LRV1}). Write $\mathcal T$ for the Bruhat--Tits tree of $\PGL_2(\Q_p)$ (see \cite[Ch. II, \S 1]{Se}), whose set of vertices (respectively, oriented edges) will be denoted by $\mathcal V$ (respectively, $\mathcal E$). The groups $B^\times$ and $(B\otimes\Q_p)^\times$ act on the left on $\mathcal T$ via $i_p$. Let $v_*\in\mathcal V$ be the distinguished vertex corresponding to the maximal order $\M_2(\Z_p)$ of $\M_2(\Q_p)$, and say that a vertex $v\in\mathcal V$ is \emph{even} (respectively, \emph{odd}) if its distance from $v_*$ is even (respectively, odd). Moreover, denote by $\hat v_*$ the vertex corresponding to $\smallmat{\Z_p}{p^{-1}\Z_p}{p\Z_p}{\Z_p}$, fix the edge $e_*=(v_*,\hat v_*)\in\mathcal E$ and say that an edge $e=(v_1,v_2)\in\mathcal E$ is \emph{even} (respectively, \emph{odd}) if $v_1$ is even (respectively, odd). Write $\mathcal E^+$ for the set of even vertices of $\mathcal T$. Finally, if $e=(v_1,v_2)$ write $\bar e$ for the reversed edge $(v_2,v_1)$.

Choose a harmonic (in the sense of \cite[Definition 4.6]{LRV1}) system of representatives $\mathcal Y_{\rm har}$ for the cosets in $\Gamma_0^D(Mp)\backslash\Gamma$, which can be written as $\mathcal Y_{\rm har}=\{\gamma_e\}_{e\in\mathcal E^+}$ with $\gamma_e\in\Gamma$ such that $\gamma_e(e)=e_*$. Set $H_E:=H_1\bigl(E(\C),\Z\bigr)$ and write $H_E^\pm$ for the $\pm$-eigenspace for complex conjugation acting on $H_E$. Fix an isomorphism $H_E^\pm\simeq\Z$. Keeping in mind that $H_E^\pm$ is (canonically isomorphic to) a quotient of $\Gamma_0^D(Mp)$, define, as in \cite[Definition 4.2]{LRV1}, the \emph{universal $1$-cochain} $\mu_{\rm univ}^\mathcal Y:\Gamma\rightarrow\mathcal F_{\rm har}(\Z)$ associated with $\mathcal Y$ by the following rules:
\begin{itemize} 
\item if $\gamma\in\Gamma$ and $e\in\mathcal E^+$ let $g_{\gamma,e}\in\Gamma_0^D(Mp)$ be defined by the equation $\gamma_e\gamma=g_{\gamma,e}\gamma_{\gamma^{-1}(e)}$, then set $\mu_{\rm univ}^\mathcal Y(\gamma)(e):=[g_{\gamma,e}]\in H_E^\pm\simeq\Z$; 
\item if $\gamma\in\Gamma$ and $e\not\in\mathcal E^+$ then set $\mu_{\rm univ}^\mathcal Y(\gamma)(e):=-\mu_{\rm univ}^\mathcal Y(\gamma)(\bar e)$. 
\end{itemize} 
Harmonic systems always exist, by \cite[Proposition 4.8]{LRV1}. The canonical generator $\boldsymbol\mu_f^\pm$ of $H^1\bigl(\Gamma,\mathcal M_0(\Q)\bigr)$ is obtained by fixing a prime number $r\nmid MDp$, applying the Hecke operator $t_r:=T_r-(r+1)$ to the class of $\mu_{\rm univ}^\mathcal Y$ in $H^1\bigl(\Gamma,\mathcal F_{\rm har}(\Z)\bigr)$ and then using the isomorphism 
\[ H^1\bigl(\Gamma,\mathcal M_0(\Q)\bigr)\simeq H^1\bigl(\Gamma,\mathcal M_0(\Z)\bigr)\otimes_\Z\Q. \]
As shown in \cite[Lemma 4.11]{LRV1}, the resulting class $\boldsymbol\mu_f^\pm$ independent of the choice of $\mathcal Y_{\rm har}$ as above. We also fix a representative $\mu_f^\pm$ of $\boldsymbol{\mu}_f^\pm$.  
  
Let $T:=\mathbb G_m$ denote the multiplicative group (viewed as a functor on commutative $\Q$-algebras). As explained in \cite[Section 6]{Gr}, the measure $\boldsymbol\mu_f^\pm$ can be used to define a lattice $L$ in $T(\C_p)$ as the image of the composition
\begin{equation} \label{int-pairing}
H_2(\Gamma,\Z)\longrightarrow H_1\bigl(\Gamma,\Div^0(\mathcal H_p)\bigr)\longrightarrow T(\C_p).
\end{equation}
Here the first map is extracted from the long exact sequence in homology associated with the short exact sequence 
\[ 0\longrightarrow\Div^0(\mathcal H_p)\longrightarrow\Div(\mathcal H_p)\xrightarrow{\rm deg}\Z\longrightarrow0 \] 
and the second is the integration map described in \cite[\S 5.1]{LRV1}. In fact, the lattice $L$ is contained in $T(\Q_p)$ and is Hecke-stable (\cite[Proposition 6.1]{LRV1}). Fix $z\in K_p-\Q_p$ and let $\boldsymbol{d}_z\in H^2\bigl(\Gamma,T(\mathbb C_p)\bigr)$ be the cohomology class represented by the $2$-cocycle 
\begin{equation} \label{d-cocycle-eq}
d_z:\Gamma\times\Gamma\longrightarrow T(K_p),\qquad(\gamma_1,\gamma_2)\longmapsto\mint_{\PP^1(\Q_p)}\frac{s-\gamma_1^{-1}(z)}{s-z}d\mu_{f,\gamma_2}^\pm(s). 
\end{equation}
(See \cite[\S 5.1]{LRV1} for the definition of the multiplicative integral.) The class $\boldsymbol{d}_z$ does not depend on the choice of the representative $\mu_f^\pm$ of $\boldsymbol\mu_f^\pm$. Write $\bar d_z$ for the composition of $d_z$ with the canonical projection onto $T(K_p)/L$. By construction, $L$ is the smallest subgroup of $T(\Q_p)$ such that $\boldsymbol d_z$ becomes trivial in $H^2\bigl(\Gamma,T(\C_p)/L\bigr)$, so there exists 
\[ \beta_z:\Gamma\longrightarrow T/L \] 
such that 
\[ \bar d_z(\gamma_1,\gamma_2)=\beta_z(\gamma_1\gamma_2)\cdot\beta_z(\gamma_1)^{-1}\cdot\beta_z(\gamma_2)^{-1} \]
for all $\g_1,\g_2\in\Gamma$ (the action of $\Gamma$ on $H_E^\pm$ used in \cite{LRV1} is the trivial one). Note that $\beta_z$ is well defined only up to elements in $\Hom(\Gamma,T/L)$. To deal with this ambiguity, recall that $\Gamma^{\rm ab}$ is a finite group (\cite[Proposition 2.1]{LRV2}). Thus, if $t$ is the exponent of $\Gamma^{\rm ab}$ then $t\cdot\beta_z$ is well defined (i.e., it depends only on the choice of a representative for $\boldsymbol\mu_f^\pm$). 

Let $H_K^+$ denote the narrow Hilbert class field of $K$ and set $G^+_K:=\Gal(H_K^+/K)$. Write $\Pic^+(\cO_K)$ for the narrow class group of the ring of integers 
$\cO_K$ of $K$. The reciprocity map of global class field theory induces an isomorphism
\begin{equation} \label{reciprocity-eq} 
\text{rec}:\Pic^+(\cO_K)\overset\simeq\longrightarrow G_K^+. 
\end{equation}
Let $\psi:K\hookrightarrow B$ be an optimal embedding of $\cO_K$ into $R_0^D(M)$, i.e. an embedding of $K$ into $B$ such that $\psi(\cO_K)=\psi(K)\cap R_0^D(M)$, 
and write $\Emb\bigl(\cO_K,R_0^D(M)\bigr)$ for the set of such embeddings. The group $\Gamma_0^D(M)$ acts on $\Emb\bigl(\cO_K,R_0^D(M)\bigr)$ by conjugation.

Recall the points $z_\psi$ and $\bar z_\psi$ defined in \S\ref{sec4.1}. By Dirichlet's unit theorem, the abelian group of units in $\mathcal O^\times_K$ of norm $1$ is free of rank one. Choose a generator $\varepsilon_K$ of this group such that $\varepsilon_K>1$ with respect to the fixed embedding $K\hookrightarrow\R$ of \S \ref{sec2.1} and set $\gamma_\psi:=\psi(\varepsilon_K)\in\Gamma_0^D(M)$.

\begin{definition} \label{darmon-def}
The \emph{Darmon points} on $T(K_p)/L$ are the points 
\[ \mathcal P_\psi:=t\cdot\beta_{z_\psi}(\gamma_\psi)\in T(K_p)/L \]
where $\psi$ varies in $\Emb\bigl(\cO_K,R_0^D(M)\bigr)$.
\end{definition}

Let $\psi$ be as before. While, thanks to the multiplicative factor $t$, the point $\mathcal P_\psi$ does not depend on the choice of a map $\beta_{z_\psi}$ splitting the $2$-cocycle $d_{z_\psi}$ modulo $L$, it might \emph{a priori} depend on the choice of the representative $\mu_f^\pm$ of $\boldsymbol\mu_f^\pm$. In fact, this is not the case: a direct calculation explained in \cite{LRV2} shows that the point $\mathcal P_\psi$ is independent of the choice of a representative of $\boldsymbol\mu_f^\pm$ and depends only on the $\Gamma_0^D(M)$-conjugacy class of $\psi$. Although, in light of this result, the symbol $\mathcal P_{[\psi]}$ would be more appropriate, in order not to burden our notation we will continue to write $\mathcal P_\psi$ for the Darmon points. However, the reader should always bear in mind that $\mathcal P_\psi=\mathcal P_{\psi'}$ whenever $\psi$ and $\psi'$ are $\Gamma_0^D(M)$-conjugate.

\subsection{Galois action on optimal embeddings} \label{Galois-emb-subsec}

As in \cite[\S 4.2]{LRV2}, for every prime $\ell$ dividing $MD$ fix orientations of $R_0^D(M)$ and $\cO_K$ at $\ell$, i.e., ring homomorphisms 
\[ \mathfrak O_\ell:R_0^D(M)\longrightarrow\mathbb F_{\ell^\delta},\qquad\mathfrak o_\ell:\mathcal O_K\longrightarrow\mathbb F_{\ell^\delta} \]
where $\delta=1$ if $\ell|M$ and $\delta=2$ if $\ell|D$. There are exactly two orientations of $\cO_K$ at every $\ell|MD$.

\begin{definition} \label{emb-def}
1) Two embeddings $\psi,\psi'\in\Emb\bigl(\cO_K,R_0^D(M)\bigr)$ are said to have \emph{the same orientation} at a prime $\ell|MD$ if $\mathfrak O_\ell\circ(\psi|_{\cO_K})=\mathfrak O_\ell\circ(\psi'|_{\cO_K})$ and are said to have \emph{opposite orientations} at $\ell$ otherwise. 

2) An embedding $\psi\in\Emb\bigl(\cO_K,R_0^D(M)\bigr)$ is said to be \emph{oriented} if $\mathfrak O_\ell\circ(\psi|_{\cO_K})=\mathfrak o_\ell$ for all primes $\ell|MD$.
\end{definition}

We denote the set of oriented optimal embeddings of $\cO_K$ into $R_0^D(M)$ by $\E\bigl(\cO_K,R_0^D(M)\bigr)$. Note that the action of $\Gamma^D_0(M)$ on $\Emb\bigl(\cO_K,R_0^D(M)\bigr)$ by conjugation restricts to an action on 
$\E\bigl(\cO_K,R_0^D(M)\bigr)$. Moreover, if $\psi\in\E\bigl(\cO_K,R_0^D(M)\bigr)$ then
\[ \psi^*:=\omega_\infty\psi\omega_\infty^{-1} \]
is in $\E\bigl(\cO_K,R_0^D(M)\bigr)$ too.

The next lemma, the proof of which we omit, collects some basic properties of (oriented) optimal embeddings.

\begin{lemma} \label{emb-lemma}
Take $\psi,\psi'\in\Emb\bigl(\cO_K,R_0^D(M)\bigr)$.
\begin{enumerate}
\item If $\psi$ and $\psi'$ have opposite orientations at a prime dividing $MD$ then at most one of them is oriented in the sense of Definition \ref{emb-def}.
\item If $\ell$ is a prime dividing $MD$ then $\psi$ and $\omega_\ell\psi\omega_\ell^{-1}$ have opposite orientations at $\ell$ and the same orientation at all the other primes dividing $MD$.
\item The embeddings $\psi$ and $\psi\circ\tau$ have opposite orientations at all the primes dividing $MD$.
\end{enumerate}
\end{lemma}
 
The principal ideal $(\sqrt{d_K})$ is a proper $\cO_K$-ideal of $K$, so we can consider its class 
$\mathfrak D_K$ in $\Pic^+(\cO_K)$; define $\sigma_K:=\text{rec}(\mathfrak D_K)\in G_K^+$. From here to the end of this subsection, for notational convenience set $R:=R_0^D(M)$. We want to recall the actions of $\Pic^+(\cO_K)$ and $G_K^+$ on the sets of $\Gamma_0^D(M)$-conjugacy classes of optimal and oriented optimal embeddings of $\cO_K$ into $R$. For details, 
see \cite[Proposition 4.2]{LRV2}.  
Let $\fr a\subset\cO_K$ be an ideal representing a class $[\fr a]\in\Pic^+(\cO_K)$ and let $\psi\in\Emb(\cO_K,R)$. Since the quaternion algebra $B$ is indefinite, the left $R$-ideal $R\psi(\fr a)$ is principal; let $a\in R$ be a generator of this ideal with positive reduced norm, which is unique up to elements in $\Gamma_0^D(M)$. The right action of $\psi(\cO_K)$ on 
$R\psi(\fr a)$ shows that $\psi(\cO_K)$ is contained in the right order of $R\psi(\fr a)$, which is equal to $a^{-1}Ra$. It can be checked that this recipe induces a well-defined action of $\Pic^+(\cO_K)$ on conjugacy classes of embeddings given by
\[ [\fr a]\star[\psi]:=\bigl[a\psi a^{-1}\bigr]\in\Emb(\cO_K,R)/\Gamma_0^D(M). \]
Notice that if $\fr a=(\sqrt{d_K})$ then we can take $a=\omega_\infty\cdot\psi(\sqrt{d_K})$, whence 
\begin{equation} \label{d-K-eq}
\fr D_K\star[\psi]=\bigl[\omega_\infty\psi\omega_\infty^{-1}\bigr]=[\psi^\ast]. 
\end{equation} 
Recalling the reciprocity map of \eqref{reciprocity-eq}, for all $\sigma\in G_c^+$ and $[\psi]\in\Emb(\cO_K,R)/\Gamma_0^D(M)$ define
\[ \sigma\star[\psi]:=\text{rec}^{-1}(\sigma)\star[\psi]. \]
In particular, it follows from \eqref{d-K-eq} that $\sigma_K\star[\psi]=[\psi^*]$ for all $\psi\in\Emb(\cO_K,R)$. Finally, for every $\sigma\in G_c^+$ and every $\psi\in\Emb(\cO_K,R)$ choose an element 
\[ \sigma\star\psi\in\sigma\star[\psi]. \]
If now $\psi$ is an oriented optimal embedding (with respect to the orientations of $\cO_K$ and $R$ fixed above) then, with notation as before, the Eichler order $a^{-1}Ra$ inherits an orientation from the one of $R$. Actually, it can be checked that we get an induced action of $\Pic^+(\cO_K)$ (and $G_K^+$) on the set $\E(\cO_K,R)/\Gamma_0^D(M)$, and this action turns out to be free and transitive. In light of this, it is convenient to describe a (non-canonical) bijection between $\E(\cO_K,R)/\Gamma_0^D(M)$ and $G_K^+$. To this end, fix once and for all an auxiliary embedding $\psi_0\in\E(\cO_K,R)$ and define $E:G_K^+\rightarrow\E(\cO_K,R)/\Gamma_0^D(M)$ by $E(\sigma):=\sigma\star[\psi_0]$. Finally, set
\[ G:=E^{-1}:\E(\cO_K,R)/\Gamma_0^D(M)\longrightarrow G_K^+, \]
which is a bijection satisfying
\begin{equation} \label{G-eq}
G([\psi^*])=\sigma_K\cdot G([\psi])
\end{equation}
for all $\psi\in\E(\cO_K,R)$.    
Now choose for every $\sigma\in G_K^+$ an embedding $\psi_\sigma\in E(\sigma)$, so that the family $\{\psi_\sigma\}_{\sigma\in G_K^+}$ is a set of representatives of the $\Gamma_0^D(M)$-conjugacy classes of oriented optimal embeddings of $\cO_K$ into $R$. If $\gamma,\gamma'\in R$ write $\gamma\sim\gamma'$ to indicate that $\gamma$ and $\gamma'$ are in the same $\Gamma_0^D(M)$-conjugacy class, and adopt a similar notation for (oriented) optimal embeddings of $\cO_K$ into $R$. For all $\sigma,\sigma'\in G_K^+$ one has
\[ \sigma\star[\psi_{\sigma'}]=\sigma\star E(\sigma')=\sigma\star\bigl(\sigma'\star[\psi_0]\bigr)=(\sigma\sigma')\star[\psi_0]=E(\sigma\sigma'), \]
which means that
\[ \sigma\star\psi_{\sigma'}\sim\psi_{\sigma\sigma'}. \]
Furthermore, since $G([\psi_\sigma^*])=\sigma\cdot G([\psi_\sigma])=\sigma_K\sigma$, by equality \eqref{G-eq}, we deduce that
\[ \psi_\sigma^*\sim \psi_{\sigma_K\sigma} \]
for all $\sigma\in G_K^+$.

To conclude this subsection, we prove a lemma which will be useful later on. As above, for every prime $\ell$ dividing $MD$ fix an orientation $\fr O_\ell:R\longrightarrow\mathbb F_{\ell^\delta}$ of the Eichler order $R$ at $\ell$.

\begin{lemma} \label{galois-transitivity-lemma}
Let $\psi,\psi'\in\Emb(\cO_K,R)$. If $\psi$ and $\psi'$ have the same orientation at all primes dividing $MD$ then there exists a unique $\sigma\in G_K^+$ such that
$\psi\sim\sigma\star\psi'$. 
\end{lemma}

\begin{proof} For every prime $\ell|MD$ choose the unique orientation of $\cO_c$ at $\ell$ such that $\psi$ and $\psi'$ become \emph{oriented} optimal embeddings. Since the action of $G_K^+$ on $\Gamma_0^D(M)$-conjugacy classes of oriented optimal embeddings is free and transitive, there exists a unique $\sigma\in G_K^+$ such that $[\psi]=\sigma\star[\psi']$, and the claim follows. \end{proof}

\begin{remark}
The element $\sigma\in G_K^+$ whose existence is established in Lemma \ref{galois-transitivity-lemma} depends on the chosen orientations $\fr O_\ell$ of $R$.
\end{remark}

\subsection{Real conjugation on Darmon points} \label{Galois-subsec}  

In the sequel we will need to understand the action of the local Galois group $\Gal(K_p/\Q_p)$ on Darmon points. As in \S \ref{sec4.1}, denote by $\tau$ the generator of $\Gal(K_p/\Q_p)$. We start with the following    

\begin{lemma} \label{lemma-galois}
$\mathcal P_{\psi\circ\tau}=\tau(\mathcal P_\psi)^{-1}$.
\end{lemma}

\begin{proof} First recall that $\bar z_\psi$ satisfies the condition $\psi(\alpha)\binom{\bar z_\psi}{1}=\bar\alpha\binom{\bar z_\psi}{1}$, so $z_{\psi\circ\tau}=\bar z_\psi$. Since $(\psi\circ\tau)(u)=\psi(u^{-1})=\gamma_\psi^{-1}$, it follows that $\mathcal P_{\psi\circ\tau}=t\cdot\beta_{\bar z_\psi}(\gamma_\psi^{-1})$. Observe that the lattice $L$, being contained in $T(\Q_p)$, inherits the trivial action of $\tau$. With notation as before, if $\beta_z$ splits $\bar d_z$ then $\tau\circ\beta_z$ splits $\bar d_{\tau(z)}$, so it follows that
\begin{equation} \label{eq++}
\mathcal P_{\psi\circ\tau}=t\cdot\beta_{\bar z_\psi}(\gamma_\psi^{-1})=t\cdot\tau\bigl(\beta_{z_\psi}(\gamma_\psi^{-1})\bigr)=\tau\bigl(t\cdot\beta_{z_\psi}(\gamma_\psi^{-1})\bigr).
\end{equation}
It remains to show that $\beta_{z_\psi}(\gamma_\psi^{-1})=\beta_{z_\psi}(\gamma_\psi)^{-1}$. To prove this equality, note that, since $\gamma_\psi(z_\psi)=z_\psi$ and $\mu_{f,\gamma}^\pm$ has total mass $0$ for all $\gamma\in\Gamma$, one has
\[ d_{z_\psi}(\gamma_\psi,\gamma)=\mint_{\PP^1(\Q_p)}\frac{s-\gamma_\psi^{-1}(z_\psi)}{s-z_\psi}d\mu_{f,\gamma}^\pm=1 \] 
for all $\gamma\in\Gamma$. Thus we obtain that
\begin{equation} \label{eq+}
1=d_{z_\psi}(\gamma_\psi,\gamma_\psi^{-1})=\beta_{z_\psi}(1)\cdot\beta_{z_\psi}(\gamma_\psi)^{-1}\cdot\beta_{z_\psi}\bigl(\gamma_\psi^{-1}\bigr)^{-1}.
\end{equation}
On the other hand, one has 
\[ 1=d_{z_\psi}(\gamma_\psi,1)=\beta_{z_\psi}(\gamma_\psi)\cdot\beta_{z_\psi}(\gamma_\psi)^{-1}\cdot\beta_{z_\psi}(1)^{-1}, \]
from which we also get $\beta_{z_\psi}(1)=1$. Substituting in \eqref{eq+} gives $\beta_{z_\psi}(\gamma_\psi^{-1})=\beta_{z_\psi}(\gamma_\psi)^{-1}$, as claimed. Equality \eqref{eq++} then becomes 
\[ \mathcal P_{\psi\circ\tau}=\tau\bigl(t\cdot\beta_{z_\psi}(\gamma_\psi)^{-1}\bigr)=\tau(\mathcal P_\psi)^{-1}, \] 
completing the proof. \end{proof}

By part (3) of Lemma \ref{emb-lemma}, if $\psi\in\Emb\bigl(\cO_K,R_0^D(M)\bigr)$ then $\psi$ and $\psi\circ\tau$ have opposite orientations at every prime dividing $MD$. Since conjugation by $\omega_{MD}$ reverses the orientation at every prime (by (2) of the same lemma) and $\omega_{MD}$ normalizes $\Gamma_0^D(M)$, Lemma \ref{galois-transitivity-lemma} ensures that there exists a unique $\sigma_\psi\in G_K^+$ such that 
\begin{equation} \label{W}
\psi\circ\tau\sim\omega_{MD}(\sigma_\psi\star\psi)\omega_{MD}^{-1}.
\end{equation}

Now we can prove

\begin{proposition} \label{prop2.9}
Let $w_{MD}$ denote the eigenvalue of $W_{MD}$ on $H_E^\pm$. Then 
\[ \tau(\mathcal P_\psi)=\mathcal P_{\sigma\star\psi}^{-w_{MD}} \] 
where $\sigma=\sigma_\psi$ is the element of $G_K^+$ appearing in \eqref{W}.
\end{proposition}

\begin{proof} We first notice that the integration map appearing in \eqref{int-pairing} is equivariant for the action of $W_{MD}$. This can be shown as in \cite[Proposition 5.1]{LRV1}, where the same result is proved for $W_p$ and $W_\infty$, so we just sketch the arguments in this case. Oserve that $\om_{MD}\in R_0^D(Mp)$ satisfyes $\G\om_{MD}=\om_{MD}\G$ and $\Gamma_0^D(Mp)\om_{MD}=\om_{MD}\Gamma_0^D(Mp)$. Let $\mathcal Y_{\rm rad}$ be a fixed radial system of representatives for $\Gamma_0^D(Mp)\backslash\Gamma$ used to compute a representative $\mu_f^\pm$ of $\boldsymbol\mu_f^\pm$ (so $\mu_{f,\g}^\pm=t_r\mu_{{\rm univ},\gamma}^{\mathcal Y_{\rm rad}}$). Introduce the system
\[ \mathcal Y'_{\rm rad}:=\bigl\{\om_{MD}^{-1}\g_{\om(e)}\om_{MD}\bigr\}_{e\in\E^+}. \]
As noticed in \cite[\S 5.2]{LRV1}, $\mathcal Y'_{\rm rad}$ is again radial and thus can be used to define a representative of $\boldsymbol\mu_f^\pm$. Since the Hecke action is defined at the cochain level, one easily verifies that 
\begin{equation} \label{**}
\mint_{\PP^1(\Q_p)}\frac{t-\omega_{MD}\gamma_1^{-1}(z_\psi)}{t-\omega_{MD}z_\psi}dt_r\mu_{{\rm univ},\omega_{MD}\gamma_2\omega_{MD}^{-1}}^{\mathcal Y_{\rm rad}}=W_{MD}\cdot\mint_{\PP^1(\Q_p)}\frac{t-\gamma_1^{-1}(z_\psi)}{t-z_\psi}dt_r\mu_{{\rm univ},\gamma_2}^{\mathcal Y'_{\rm rad}}
\end{equation} 
for all $\gamma_1,\gamma_2\in\Gamma$. To simplify notations, set $\mu'_{\g}:=t_r\mu_{{\rm univ},\gamma}^{\mathcal Y'_{\rm rad}}$ and $\psi':=\omega_{MD}(\sigma\star\psi)\omega_{MD}^{-1}$, with $\sigma=\sigma_\psi$ as in \eqref{W}. Then $z_{\psi'}=\omega_{MD}(z_\psi)$ and $\bar z_{\psi'}=\omega_{MD}(\bar z_\psi)$. By \eqref{**}, 
if $\beta'_{z_\psi}$ splits the $2$-cocycle
\[ (\gamma_1,\gamma_2)\longmapsto\mint_{\PP^1(\Q_p)}\frac{t-\gamma_1^{-1}(z_\psi)}{t-z_\psi}d\mu'_{\g_2} \]
then $\gamma\mapsto w_{MD}\beta'_{z_\psi}(\omega_{MD}^{-1}\gamma\omega_{MD})$ splits the $2$-cocycle $d_{\omega_{MD}z_\psi}=d_{z_{\psi'}}$ defined in \eqref{d-cocycle-eq} (here use the fact that $W_{MD}$ acts on $T $ as multiplication by $w_{MD}$). If follows that 
\[ w_{MD}\cdot\beta'_{z_\psi}(\omega_{MD}^{-1}\gamma\omega_{MD})=\beta_{z_{\psi'}}(\g)+\varphi(\g) \] 
where $\varphi:\Gamma\rightarrow T(K_p)/L_T $ is a homomorphism. Therefore we get 
\begin{equation} \label{eq1+}
\mathcal P_{\psi'}=t\cdot\beta_{z_{\psi'}}(\gamma_{\psi'})=w_{MD}t\cdot\beta'_{z_\psi}(\omega_{MD}^{-1}\gamma_{\psi'}\omega_{MD})=w_{MD}t\cdot\beta'_{z_\psi}(\gamma_{\sigma\star\psi}).
\end{equation}
As remarked in \S \ref{sec2.5}, $t\cdot\beta'_{z_\psi}(\gamma_{\sigma\star\psi})=t\cdot\beta_{z_\psi}(\gamma_{\sigma\star\psi})$, hence \eqref{eq1+} implies that 
\begin{equation} \label{eq1++} 
\mathcal P_{\psi'}=\mathcal P_{\sigma\star\psi}^{w_{MD}}.
\end{equation} 
From Lemma \ref{lemma-galois} and \eqref{W} we obtain
\begin{equation} \label{eq2+}
\tau(\mathcal P_\psi)=\mathcal P_{\psi'}^{-1},
\end{equation}
and then combining \eqref{eq1++} with \eqref{eq2+} yields the result. \end{proof}

\subsection{Darmon points on elliptic curves} \label{sec2.2}

Define a $p$-adic logarithm on the quotient $T(\C_p)/L$ as follows. As in the introduction, let $q\in p\Z_p$ be Tate's period for $E$ at $p$ and let $\log_q$ be the branch of the $p$-adic logarithm such that $\log_q(q)=0$. Since $T(\C_p)\simeq\C_p^\times$, the fact that $\log_q$ induces a map on $T(\C_p)/L$ is an immediate consequence of the following

\begin{lemma} \label{log-lemma}
$\log_q(L)=0$.
\end{lemma}

\begin{proof} By \cite[Theorem 7.16]{LRV1}, the lattices $L$ and $\langle q\rangle$ are homothetic, i.e.,
$L\cap\langle q\rangle$ has finite index both in $L$ and in $\langle q\rangle$. 
If $n:=[L:L\cap\langle q\rangle]$ then $L^n\subset\langle q\rangle$, so $\log_q(L^n)=0$ by the choice of 
$\log_q$. Since $\log_q$ is a homomorphism from $\C_p^\times$ to $\C_p$, the lemma is proved. \end{proof}

From now on, by an abuse of notation we will use the symbol $\log_q$ for both maps induced by $\log_q$ on $T(\bar K_p)/L$ 
and $\bar K_p^\times/\langle q\rangle$; these two maps are $K_p$-valued on the $K_p$-rational points. 

As in the proof of the lemma above, set $n:=[L:L\cap\langle q\rangle]$. 
Raising to the $n$-th power induces a Galois-equivariant isogeny of $p$-adic tori
\begin{equation} \label{tori-isogeny-eq}
T(\bar K_p)/L\longrightarrow\bar K_p^\times/\langle q\rangle 
\end{equation}
defined over $K_p$, and then composing \eqref{tori-isogeny-eq} with Tate's analytic uniformization $\Phi_{\mathrm{Tate}}$ introduced in \eqref{Tate} yields a Galois-equivariant isogeny
\[ \varphi_E:T(\bar K_p)/L\longrightarrow E(\bar K_p). \]
As in \cite[p. 346]{BD-annals}, define a homomorphism 
\[ \log_E:E(K_p)\longrightarrow K_p \]
by the rule $\log_E(P):=\log_q\bigl(\Phi_{\mathrm{Tate}}^{-1}(P)\bigr)$ for all $P\in E(K_p)$. It follows that 
\begin{equation} \label{log-comparison-eq}
\log_E\bigl(\varphi_{E}(x)\bigr)=n\log_q(x)\qquad\text{for all $x\in T(K_p)/L$}. 
\end{equation}
For every $\psi\in\Emb\bigl(\cO_c,R_0^D(M)\bigr)$ define
\[ P_\psi^\pm:=\varphi_E\bigl(\mathcal P_\psi\bigr)\in E(K_p). \]

\begin{remark}
Since the same is true of $\mathcal P_\psi$, the point $P_\psi^\pm$ is independent of the representative of $\boldsymbol\mu_f^\pm$, and $P_\psi^\pm=P_{\psi'}^\pm$ if $\psi$ and $\psi'$ are $\Gamma_0^D(M)$-conjugate. However, as made explicit by the notation, $P_\psi^\pm$ \emph{does} depend on the choice of sign $\pm$.
\end{remark}

\section{Explicit formulas for Darmon points} \label{sec-control} 

The aim of this section is to use a control theorem of Greenberg--Stevens type in the quaternionic setting (Theorem \ref{control-thm}) to obtain a canonical lift $\boldsymbol{\tilde\mu}_f^\pm$ of $\boldsymbol{\mu}_f^\pm$. The role of this explicit lift is twofold: on the one hand, it is used to perform computations with Darmon points (cf. Corollary \ref{coro}); on the other hand, in Section \ref{L-function} it allows us to define a $p$-adic $L$-function which, thanks to Popa's work (\cite{Po}), encodes (for suitable integers $k$) special values of the $L$-functions of the Hida forms $f_k$ over $K$ twisted by genus characters of $K$. 

\subsection{Modular forms on quaternion algebras} \label{sec2.4} 

For any ring $R$ and any integer $n\geq0$ let $P_n(R):={\rm Sym}^n(R)$ denote the $R$-module of homogeneous polynomials $P(x,y)$ in two variables $x,y$ of degree $n\geq0$ with coefficients in $R$. It is equipped with a right action of the group $\GL_2(R)$ given by the rule  
\[ (P|\gamma)(x,y):=P(ax+by,cx+dy)\qquad\text{for $\gamma=\smallmat abcd$}. \]
The $R$-linear dual $V_n(R)$ of $P_n(R)$ is then endowed with a left action of $\GL_2(R)$ by the formula
\[ (\gamma\cdot\phi)(P):=\phi(P|\gamma). \]
Finally, if $R$ is a field extension of $F$ then $P_n(R)$ (respectively, $V_n(R)$) is equipped with a right (respectively, left) action of $B^\times$ via the map $i_F$ and the rules above.

Let $g$ be a cusp form of level $\Gamma_0^D(Mp)$ and weight $k$, suppose that $g$ is an eigenform for all the Hecke operators in the usual Hecke algebra $\fr h_k$ and write $\lambda:\mathfrak h_k\rightarrow\C$ for the corresponding morphism. Let $L$ be a subfield of $\C$ containing the image of $\lambda$. If $\iota$ denotes complex conjugation, define $H^1\bigl(X_0^D(Mp),V_{k-2}(L)\bigr)^{g,\pm}$ to be the $L$-subspace of $H^1\bigl(X_0^D(Mp),V_{k-2}(L)\bigr)$ consisting of the elements $\xi$ such that $\xi|T=\lambda(T)\xi$ for all $T\in\mathfrak h_{k}$ and $\xi|\iota=\pm\xi$. By \cite{MS}, we know that 
\[ \dim_L\Big(H^1\bigl(X_0^D(Mp),V_{k-2}(L)\bigr)^{g,\pm}\Big)=1. \]
Recall that there is a canonical isomorphism 
\[ H^1\bigl(X_0^D(Mp),V_{k-2}(\C)\bigr)\simeq H^1\bigl(\Gamma_0^D(Mp),V_{k-2}(\C)\bigr) \] 
which is equivariant for the action of the involution $\iota$. Then the complex vector space $H^1\bigl(\Gamma_0^D(Mp),V_{k-2}(\C)\bigr)^{g,\pm}$ is spanned by the projection on the $\pm$-eigenspace for $\iota$ of the class represented by the $1$-cocycle $\gamma\mapsto v(g)_\gamma$ with 
\[ v(g)_\gamma(P):=\int_\tau^{\gamma(\tau)}g(z)P(z,1)dz \]
(this class does not depend on the choice of $\tau\in\mathcal H$). For details, see \cite[\S 8.2]{Sh}. 

\subsection{Hida families} 

Recall the group $\mathcal X$ in \eqref{X} and the Hida family $f_\infty=\sum_{n=1}^\infty a_n(\kappa)q^n$ of $f$ in
\eqref{q-exp}, 
where the $p$-adic analytic functions $a_n$ are all defined in a $p$-adic neighbourhood $\U_f$ of $2$ in $\mathcal X$. For a fixed even integer $k\geq2$ in $\U_f$ write $f_k$ for the normalized eigenform of weight $k$ on $\Gamma_0(MDp)$ whose $q$-expansion is obtained by setting $\kappa=k$ in $f_\infty$; in particular, $f_2=f$ (see \cite[\S 1.2]{BD} for details). Shrinking it if necessary, we can assume that $\U_f$ is contained in the residue class of $2$ modulo $p-1$. For any even integer $k\in\U_f$ denote by $F_k$ the finite extension of $\Q_p$ generated by the Fourier coefficients of $f_k$. In the notations of \S \ref{sec2.4}, define the one-dimensional $F_k$-vector space 
\[ \W_{f_k}^\pm:=H^1\bigl(\Gamma_0^D(Mp),V_{k-2}(F_k)\bigr)^{f_k,\pm}, \] 
then fix a generator $\boldsymbol\phi^\pm_k$ of $\W_{f_k}^\pm$ as well as a representative $\phi_k^\pm$ of $\boldsymbol\phi_k^\pm$.

\subsection{The control theorem} \label{subsec-control} \label{sec-lifting} 

Let $\X:=(\Z_p^2)^\prime$ denote the set of primitive vectors in $\Y:=\Z_p^2$. For any $\Z_p$-module $G$ let $\tilde\D(G)$ be the group of $G$-valued measures on $\Y$ and let $\D(G)$ be the subgroup consisting of those measures which are supported on $\X$. The left action of $\Sigma:=\GL_2(\Q_p)\cap\M_2(\Z_p)$ 
on $\tilde\D(G)$ considered in \cite{LRV1} is given by the integration formula 
\[ \int_\Y\varphi(x,y)d(\gamma\cdot\nu)(x,y)=\int_\Y\varphi\bigl(\gamma\cdot(x,y)\bigr)d\mu(x,y) \] 
with 
$\gamma\cdot(x,y):=(ax+by,cx+dy)$ 
for all step functions $\varphi:\Y\rightarrow G$, all measures $\nu\in\tilde\D(G)$ and all matrices $\gamma=\smallmat abcd\in\Sigma$. This action also induces an action of $\Sigma$ on $\D(G)$ (see \cite[Lemma 7.4]{LRV1}). Denote by $\pi:\X\rightarrow \PP^1(\Q_p)$ the fibration defined by $(a,b)\mapsto a/b$, so that we obtain a canonical map 
\[\pi_*:H^1\big(\Gamma_0^D(M),\D(\Q_p)\big)\longrightarrow H^1\big(\Gamma_0^D(M),\mathcal M(\Q_p)\big).\]
To simplify the notation, write $\D:=\D(\Z_p)$. The group $H^1\bigl(\Gamma_0^D(M),\D\bigr)$ comes equipped with an action of the Hecke algebra $\mathcal H(M)$, which allows us to define its \emph{ordinary submodule} as
\[ \W:=\bigcap _{n=1}^\infty H^1\bigl(\Gamma_0^D(M),\D\bigr)\big|T_p^n \] 
where $T_p\in\mathcal H(M)$ is the $p$-th Hecke operator. One can show that $T_p$ acts invertibly on $\W$ and that $H^1\bigl(\Gamma_0^D(M),\D\bigr)$ decomposes as a direct sum of $\W$ and a submodule on which $T_p$ acts nilpotently. Finally, let 
\[\Lambda:=\Z_p[\![1+p\Z_p]\!]\simeq\Z_p[\![T]\!]\] be the Iwasawa algebra. Following the constructions in \cite[Sections 1 and 5]{GS}, one endows $\W$ with a structure of (finitely generated) $\Lambda$-module. 

The following Control Theorem is probably known to experts. To state this result, let us introduce a minimal set of notations. For details, proofs, and a dictionary between classical Hida families and their quaternionic counterparts, see \cite[Sections 5 and 6]{LV}. 

As in \cite[\S 5.3]{LV}, let $\mathcal R$ be the integral closure of $\Lambda$ in the primitive component to which $f$ belongs. The ring $\mathcal R$ is a complete local noetherian domain, finitely generated as a $\Lambda$-module and equipped with a canonical map $\eta:\mathfrak h_\infty^\ord\rightarrow\mathcal R$ of $\Lambda$-algebras (where $\mathfrak h_\infty^\ord$ is Hida's ordinary Hecke algebra: see \cite[\S 5.1]{LV}). If $k\in\U_f$ is an even positive integer there exists a continuous morphism $\vartheta_k:\mathcal R\rightarrow \bar\Q_p$ whose restriction to group-like elements in $\Lambda$ is the character $x\mapsto x^{k-2}$ and such that $f_k=\sum_{n=1}^\infty(\vartheta_k\circ\eta)(T_n)q^n$, where $T_n\in\mathfrak h_\infty^\ord$ is the $n$-th Hecke operator. Set $\mathfrak p_k:=\ker(\vartheta_k)$; if $\mathcal R_{\mathfrak p_k}$ is the localization of $\mathcal R$ at $\mathfrak p_k$, we have $\mathcal R_{\mathfrak p_k}/\mathfrak p_k\mathcal R_{\mathfrak p_k}\simeq F_k$. Define 
\[ \W_{\mathfrak p_k}^\pm:=\W^\pm\otimes_\Lambda\mathcal R_{\mathfrak p_k} \] 
where $\W^\pm$ denotes the $\pm$-eigenmodule for the action of $W_\infty$ on $\W$.

Let $\boldsymbol\nu\in\W$ and fix a representative $\nu$ of $\boldsymbol\nu$ and a field extension $L$ of $\Q_p$. The function  
\[ \gamma\longmapsto\rho_k(\nu_\gamma):=\bigg(P(x,y)\mapsto\int_{\Z_p\times\Z_p^\times}P(x,y)d\nu_\gamma\bigg) \] 
defined on $\Gamma_0^D(Mp)$ with values in $V_{k-2}(L)$ is a $1$-cocycle. Its class in $H^1\bigl(\Gamma_0^D(Mp),V_{k-2}(L)\bigr)$, denoted $\rho_k(\boldsymbol\nu)$, does not depend on the choice of $\nu$.

\begin{theorem} \label{control-thm}
Let $k\in\U_f$ be an even positive integer. For any choice of sign $\pm$ the $\mathcal R_{\mathfrak p_k}$-module $\W_{\mathfrak p_k}^\pm$ is free of rank $1$ and $\rho_k$ induces an isomorphism $\rho_k:\W^\pm_{\mathfrak p_k}/\mathfrak p_k\W_{\mathfrak p_k}^\pm\overset\simeq\longrightarrow\W_{f_k}^\pm$.  
\end{theorem}

The proof of Theorem \ref{control-thm} follows the arguments in \cite[\S 5 and \S 6]{GS} closely, and the articles \cite{BL} and \cite{DG} provide generalizations of the results of \cite{GS} to the quaternionic setting we are working in. However, for lack of precise references for the statement we need, a full proof of a more general version of it will be given in a subsequent paper.

\begin{proposition} \label{prop-pi-2}
There exists $\boldsymbol{\tilde\mu}_f^\pm\in\W^\pm\otimes_{\Z_p}\Q_p$ with the following properties: 
\begin{enumerate} 
\item $\pi_*(\boldsymbol{\tilde\mu}_f^\pm)=\res^\Gamma_{\Gamma_0^D(M)}(\boldsymbol\mu_f^\pm)$ in $H^1\bigl(\Gamma_0^D(M),\mathcal M_0(\Q_p)\bigr)$;
\item at the cost of replacing $\U_f$ with a smaller neighbourhood, for every $k\in\U_f$ there exists $\lambda_B^\pm(k)\in F_k^\times$ such that $\boldsymbol\phi_k^\pm=\lambda_B^\pm(k)\rho_k(\boldsymbol{\tilde\mu}_f^\pm)$ in $\W_{f_k}^\pm$.
\end{enumerate}
\end{proposition}

\begin{proof} In light of Theorem \ref{control-thm}, we can (and do) choose $\boldsymbol{\Phi}^\pm\in\W^\pm$ such that $\rho_2(\boldsymbol{\Phi}^\pm)=\lambda_1^\pm\boldsymbol{\phi}_2^\pm$ with $\lambda_1^\pm\in\Q_p^\times$. Fix a representative $\Phi^\pm$ of $\boldsymbol\Phi^\pm$ and note that $\rho_2(\boldsymbol\Phi^\pm)$ is represented by the cocycle $\gamma\mapsto(\pi_*(\Phi_\gamma^\pm))(\Z_p)$ on $\Gamma_0^D(Mp)$. By \cite[Proposition 25]{Gr}, there is a canonical isomorphism
\[ \mathcal S:H^1\bigl(\Gamma,\mathcal M_0(\Q_p)\bigr)^{f,\pm}\overset\simeq\longrightarrow\W_f^\pm, \] 
and thus there exists $\lambda_2^\pm\in\Q_p^\times$ such that $\mathcal S(\boldsymbol{\mu}_f^\pm)=\lambda_2^\pm\boldsymbol{\phi}_2^\pm$. Moreover, observe that $\mathcal S({\boldsymbol\mu}_f^\pm)$ is represented by the cocycle $\gamma\mapsto\mu_{f,\gamma}^\pm(\Z_p)$ on $\Gamma_0^D(Mp)$. 

Let $\lambda^\pm:=\lambda_2^\pm/\lambda_1^\pm$, define $c_\gamma:=\lambda^\pm\pi_*(\Phi_\gamma^\pm)-\mu_{f,\gamma}^\pm$ for $\gamma\in\Gamma_0^D(M)$ and set
$\boldsymbol c:=[c]$ in $H^1\bigl(\Gamma_0^D(M),\mathcal M(\Q_p)\bigr)$. Using the surjectivity of $\pi_*$ proved in \cite[Proposition 7.7]{LRV1}, choose $\boldsymbol{\tilde c}\in\W^\pm\otimes_\Z\Q$ such that $\pi_*(\boldsymbol{\tilde c})=\boldsymbol{c}$. (Note that one can choose $\boldsymbol{\tilde c}$ in the 
ordinary submodule because the arrows between Tate modules at the end of the proof of \cite[Proposition 7.7]{LRV1} respect ordinary summands; the same remark applies 
to the $\pm$-eigenmodules for $W_\infty$). Since $[\gamma\mapsto c_\gamma(\Z_p)]=0$ in $\W_f^\pm$, one has $\rho_2(\boldsymbol{\tilde c})=0$. Therefore the class $\boldsymbol{\tilde\mu}_f^\pm:=\lambda^\pm\boldsymbol{\Phi}^\pm-\boldsymbol{\tilde c}$ satisfies $\rho_2(\boldsymbol{\tilde\mu}_f^\pm)=\mathcal S(\boldsymbol{\mu}_f^\pm)$ and  $\pi_*(\boldsymbol{\tilde\mu}_f^\pm)=\res^\Gamma_{\Gamma_0^D(M)}(\boldsymbol{\mu}_f^\pm)$, which proves part (1).

As for part (2), by Theorem \ref{control-thm} it is enough to observe that, since the intersection of infinitely many prime ideals of height $1$ of $\mathcal R$ is trivial, $\rho_k(\boldsymbol{\tilde\mu}_f^\pm)$ is $0$ only for finitely many integers $k$, among which the integer $2$ does not appear because $\rho_2(\boldsymbol{\tilde\mu}_f^\pm)=\mathcal S(\boldsymbol{\mu}_f^\pm)$. \end{proof}

\begin{definition} \label{defi-domain}
The open set $\U_f$ of Proposition \ref{prop-pi-2} is the \emph{domain of analiticity} of $\boldsymbol{\tilde\mu}_f^\pm$.
\end{definition}

\subsection{Explicit formulas for Darmon points} \label{sec-explicit}

As explained in \cite[\S 7.3]{LRV1}, we can choose a representative $\tilde\mu_f^\pm$ of $\boldsymbol{\tilde\mu}_f^\pm$ in Proposition \ref{prop-pi-2} so that 
$\pi_*(\tilde\mu_{f,\gamma}^\pm)=\mu_{f,\gamma}^\pm$ for all $\gamma\in\Gamma_0^D(M)$ (here the map $\pi_*$ is viewed as defined on cocycles). From now on we fix such a choice of $\tilde\mu_f^\pm$.   

Let $\log:\C_p^\times\rightarrow\C_p$ be any branch of the $p$-adic logarithm; such a branch is determined by choosing $\theta\in\C_p^\times$ with $|\theta|_p<1$ (where $|\cdot|_p$ is the $p$-adic absolute value on $\C_p$ normalized so that $|p|_p=1/p$), setting $\log(\theta):=0$ and requiring that $\log$ be a homomorphism. In particular, since $x-z_\psi y\in(\mathcal O_K\otimes\Z_p)^\times$ when $(x,y)\in\X$, the value $\log(x-z_\psi y)$ for $(x,y)\in\X$ does not depend on the choice of $\log$ (see also the discussion in the proof of \cite[Theorem 2.5]{BD-annals}). 

To simplify the notations, let $G_0:=\Gamma_0^D(Mp)$ and $G_1:=\Gamma_0^D(M)$. Define $G_2=\hat\Gamma_0^D(M):=\omega_p^{-1}G_1\omega_p$. Recall that $\Gamma$ is the amalgamated product $\Gamma=G_1*_{G_0}G_2$ with respect to the natural inclusion $G_0\subset G_1$ and the injection $\varsigma:G_0\hookrightarrow G_2$ sending $x$ to $\omega_px\omega_p^{-1}$ (cf. \cite[\S 4.1]{LRV1}).  

\begin{theorem} \label{thm2.11}
$\displaystyle{\log(\mathcal P_\psi)=-t\cdot\int_\X\log(x-z_\psi y)d{{\tilde\mu}}_{f,{\gamma_\psi}}^\pm}\pmod{\log(L)}$. 
\end{theorem} 

\begin{proof} Define the $1$-cochain $\rho_\gamma:=-\int_\X\log(x-z_\psi y)d\tilde\mu_{f,\gamma}^\pm$ for $\gamma\in G_1$ and let $\hat\rho_\gamma$ be defined from $\rho_\gamma$ for $\gamma\in G_2$ as in \cite[\S 7.5]{LRV1}. Recall that $\pi_*({\tilde\mu}_f^\pm)=\mu_f^\pm$. A combination of this equality with the calculations in \cite[\S 7.6]{LRV1} shows that $\log(\boldsymbol{d}_z)=\Delta(\boldsymbol\rho-\boldsymbol{\hat\rho})$, where $\Delta$ is the connecting map appearing in the Mayer--Vietoris exact sequence and $\boldsymbol{\rho}-\boldsymbol{\hat\rho}\in H^1(G_0,K_p)$ is represented by $\rho-\hat\rho$. We already know that, modulo $\log(L)$, the $2$-cocycle $\log(\boldsymbol{d}_z)$ is split by $\log(\beta_z)$. Hence
\[ \delta\big(\log(\beta_z)\big)=\Delta(\rho-{\hat \rho})\pmod{\log(L)}. \] 
Since $t$ annihilates $\Gamma^{\rm ab}$, the claim is then a consequence of the following lemma. \end{proof}

\begin{lemma} Let $A$ be an abelian group with trivial $\Gamma$-action and consider the portion of the Mayer--Vietoris exact sequence given by 
\[ H^1(G_1,A)\oplus H^1(G_2,A)\longrightarrow H^1(G_0,A)\overset\Delta\longrightarrow H^2(G,A)\longrightarrow H^2(G_1,A)\oplus H^2(G_2,A).\] 
Assume that 
\begin{itemize} 
\item $c=\Delta(\rho)$ for some $\rho\in H^1(G_0,A)$, so there are $1$-cochains $\theta_1\in C^1(G_1,A)$ and $\theta_2\in C^1(G_2,A)$ such that $c_{|G_1}=\delta(\theta_1)$ and $c_{|G_2}=\delta(\theta_2)$, where $\delta$ is the connecting map in cohomology and $\rho=\theta_{1|G_0}-\theta_{2|G_0}$;
\item $c$ is trivial in $H^2(G,A)$, so $c=\delta(b)$ for some $1$-cochain $b\in C^1(G,A)$.
\end{itemize}
Then $b_{|G_1}=\theta_1+\varphi_{1|G_1}$ and $b_{|G_2}=\theta_2+\varphi_{2|G_2}$ where $\varphi_1,\varphi_2:G\rightarrow A$ are group homomorphisms.
\end{lemma} 

\begin{proof} We only prove the existence of $\varphi_1$, since the result for $\varphi_2$ can be shown in the same way up to the non-trivial permutation of the set of indices $\{1,2\}$. First of all, since $c_{|G_1}=\delta(\theta_1)=\delta(b_{|G_1})$ and all the groups involved act trivially on $A$, it is immediate to deduce that
\[ b_{|G_1}=\theta_1+\psi_1 \]
where $\psi_1:G_1\rightarrow A$ is a group homomorphism. We want to prove that $\psi_1$ extends to a homomorphism $\varphi_1:G\rightarrow A$. To this end, consider the natural isomorphism $\theta:G_1\overset\simeq\rightarrow G_2$ defined by $x\mapsto\omega_px\omega_p^{-1}$ and define the homomorphism
\[ \tilde\psi_1:=\psi_1\circ\theta^{-1}:G_2\longrightarrow A \]
It follows that $\tilde\psi_1\circ\varsigma=\psi_{1|G_0}$. But then the universal property of the amalgamated product (see \cite[Ch. I, \S 1]{Se}) ensures that there exists a homomorphism $\varphi_1:G\rightarrow A$ such that $\varphi_{1|G_1}=\psi_1$, as was to be shown. \end{proof}

Now recall the notation in \S\ref{sec2.2}. 

\begin{corollary} \label{coro} 
$\displaystyle{\log_E(P_\psi)=-nt\cdot\int_\X\log_q(x-z_\psi y)d\tilde\mu_{f,\gamma_\psi}^\pm}$. 
\end{corollary} 

\begin{proof} Immediate from Theorem \ref{thm2.11} and formula \eqref{log-comparison-eq}. \end{proof}

\section{Special values of $L$-functions} \label{L-function} 

In this section we use the class $\boldsymbol{\tilde\mu}_f^\pm$ to define a $p$-adic $L$-function $L_p(f_\infty/K,\chi,\kappa)$ associated with $f_\infty$, $K$ and a genus character $\chi$ of $K$ and defined on the domain of analyticity $\U_f$ of $f$ (cf. Definition \ref{defi-domain}). Using Popa's results in \cite{Po}, we establish a link between the values of $L_p(f_\infty/K,\chi,\kappa)$ at even positive integers $\kappa=k$ and the special values of the complex $L$-functions of the forms $f_k$ over $K$ twisted by $\chi$. 

\subsection{Modular forms and $1$-cocycles} \label{sec4.1-bis} 

First, we extend the notation employed in \S\ref{sec2.4}. For any cusp form $g$ of level $\Gamma_0^D(Mp)$ and weight $k\geq2$, any point $\tau\in\mathcal H$ and any element $\gamma\in B^\times$ of positive norm define $I(g,\tau,\gamma)\in V_{k-2}(\C)$ by the formula  
\[ I(g,\tau,\gamma)\bigl(P(x,y)\bigr):=\int_\tau^{\gamma(\tau)}g(z)P(z,1)dz \] 
for all $P\in P_{k-2}(\C)$. As in \S \ref{sec2.4}, let $B^\times$ act on the left on $V_n(\C)$ via $i_F$. A direct computation shows that  
\begin{equation} \label{eq-31}
\bigl(b\cdot I(g|b,\tau,\gamma)\bigr)(P)=\det\bigl(i_F(b)\bigr)^{k-2}I\bigl(g,b(\tau),b\gamma b^{-1}\bigr)(P) 
\end{equation}
for all $g$ and $\tau$ as above, all $b,\gamma\in B^\times$ of positive norm and all $P\in P_{k-2}(\C)$. The cohomology class in $H^1\bigl(\Gamma_0^D(Mp),V_{k-2}(\C)\bigr)$ represented by the cocycle $\gamma\mapsto I(g,\tau,\gamma)$ does not depend on the choice of $\tau$. Denote $\gamma\mapsto I^\pm(g,\tau,\gamma)$ the $\pm$-component of $\gamma\mapsto I(g,\tau,\gamma)$ under the action of $\iota$; the class represented by this cocycle is again independent of the choice of $\tau$. 

If $f_k^{\rm JL}$ is a form on $\Gamma_0^D(Mp)$ associated with $f_k$ by the Jacquet--Langlands correspondence, it follows that the cohomology class represented by 
$\gamma\mapsto I^\pm(f_k^{\rm JL},\tau,\gamma)$ is a multiple of the generator $\boldsymbol{\phi}_k^\pm$ introduced in \S\ref{sec-lifting}. From now on, for a fixed sign $\pm$ we normalize the choice of $f_k^{\rm JL}=f_{k,\pm}^{\rm JL}$ in such a way that 
\begin{equation} \label{normalization}
\bigl[\gamma\mapsto I^\pm\bigl(f_k^{\rm JL},\tau,\gamma\bigr)\bigr]=\boldsymbol{\phi}_k^\pm
\end{equation} 
in $H^1\bigl(\Gamma_0^D(Mp),V_{k-2}(\C)\bigr)^{f_k,\pm}$, where we identify $\boldsymbol{\phi}_k^\pm\in\W_k^\pm$ with its image in this cohomology group via the embedding in \eqref{embeddings}.

Let $L=\Q(\sqrt{d_L})$ be a real quadratic field with discriminant $d_L$ such that 
\begin{itemize}
\item all the primes dividing $M$ are split in $L$; 
\item all the primes dividing $D$ are inert in $L$; 
\item the prime $p$ is unramified in $L$.
\end{itemize} 
Let $S$ denote either $M$ or $Mp$ and let $\psi:L\hookrightarrow B$ be an optimal embedding of the maximal order $\mathcal O_L$ of $L$ into $R_0^D(S)$ (note that such an embedding exists unless $S=Mp$ and $p$ is inert in $L$). As done in \S \ref{sec4.1} for $L=K$, we may consider the element $i_F\bigl(\psi(\sqrt{d_L})\bigr)=\smallmat abc{-a}$ and the quadratic form 
\[ Q_\psi(x,y)=cx^2-2axy-by^2. \] 
In order to uniformize our notation, we assume that $L\neq F$, so that $c\neq0$. Moreover, for every even integer $k\geq2$ define 
\[ Q^{(k)}_\psi(x,y):=Q_\psi(x,y)^\frac{k-2}{2}. \] 
As before, for $\gamma=\smallmat abcd\in\Sigma:=\GL_2(\Q_p)\cap\M_2(\Z_p)$ write $\gamma(x,y):=(ax+by,cx+dy)$. Then   
\begin{equation} \label{eq-21}
Q_{\gamma\psi\gamma^{-1}}(x,y)=\det(\gamma)\cdot Q_\psi\bigl(\gamma^{-1}(x,y)\bigr)
\end{equation} 
for all $\gamma\in\Sigma$.

Fix an embedding $L\hookrightarrow \R$ (coinciding with the one fixed in \S \ref{sec2.1} when $L=K$) and a fundamental unit $\varepsilon_L$ of $L$ such that $\varepsilon_L>1$ under this embedding, then set $\gamma_\psi:=\psi(\varepsilon_L)$. An easy computation shows that if $k\geq2$ is the weight of $g$ then 
\begin{equation} \label{pm-I-short-eq}
I(g,\psi):=I(g,\tau,\gamma_\psi)\Big(Q^{(k)}_\psi(x,y)\Big) 
\end{equation} 
does not depend on $\tau$, justifying the notation. Moreover, one checks that $I(g,\psi)$ depends only on the $\Gamma_0^D(S)$-conjugacy class of $\psi$. We define $I^\pm(g,\psi)$ in an analogous manner. 

\subsection{Preliminary calculations} 

Before proceeding, in this subsection we collect a number of preliminary calculations that will be used in the proof of the main results of this section, namely Theorems \ref{teo4.17} and \ref{thm-factor}.  

Let $\nu|\U$ denote the restriction of a measure $\nu$ on $\Y$ to the compact open subset $\U\subset\Y$. Define the compact open subsets  
\[ \X_\infty:=\Z_p^\times\times p\Z_p,\qquad\X_{\rm aff}:=\Z_p\times\Z_p^\times \]
of $\X$, so that $\X=\X_{\rm aff}\coprod\X_\infty$. 

\begin{lemma} \label{lemma4.3} 
With self-explaining notation, for every $\boldsymbol\nu\in H^1\bigl(\Gamma_0^D(M),\D\bigr)$ one has 
\[ \res(\boldsymbol\nu)|\X_\infty=\bigl(\res(\boldsymbol\nu)|\X_{\rm aff}\bigr)|U_p^{-1}W_p. \]   
\end{lemma}
\begin{proof} A variation on the proof of \cite[Lemma 7.11]{LRV1}, where a different convention for the Hecke action is adopted. The details are omitted. \end{proof}
 
\begin{notation} 
In the following, we often let $[\gamma\mapsto c_\gamma]$ or even $[c_\gamma]$ denote the class represented in $H^1(G,M)$ by a cocycle $c\in Z^1(G,M)$. Moreover, we write $c\equiv c'$ or (by abuse of notation) $c_\gamma\equiv c_\gamma'$ if the two cocycles $c$ and $c'$ are cohomologous, i.e., differ by a coboundary.
\end{notation}

\begin{lemma} \label{lemma-4} 
Let $L$ be a field, let $\boldsymbol\nu\in H^1\bigl(\Gamma_0^D(Mp),\D\bigr)$ and fix a representative $\nu$ of $\boldsymbol\nu$. Then 
\[ \bigg[P\longmapsto\int_{\Z_p^\times\times\Z_p^\times}P(x,y)d\nu_\gamma\bigg]=\bigg[P\longmapsto \int_{\X_{\rm aff}}P(x,y)d\nu_\gamma-\int_{\X_\infty}(P|\alpha)(x,y)d\nu_{\alpha^{-1}\gamma\alpha}\bigg] \] 
in $H^1\big(\Gamma_0^D(Mp),V_{k-2}(L)\big)$, where $\alpha\in\Gamma_0^D(M)$ satisfies $\alpha\X_\infty=p\Z_p\times\Z_p^\times$.  
\end{lemma}

\begin{proof} Using the decomposition $\Z_p^\times\times\Z_p^\times=\X_{\rm aff}-(p\Z_p\times\Z_p^\times)$, write
\[ \int_{\Z_p^\times\times\Z_p^\times}P(x,y)d\nu_\gamma=\int_{\X_{\rm aff}}P(x,y)d\nu_\gamma-\int_{p\Z_p\times\Z_p^\times}P(x,y)d\nu_\gamma. \]
Let $\alpha$ be as in the statement of the lemma. First of all, observe that 
\[ \int_{p\Z_p\times\Z_p^\times}P(x,y)d\nu_\gamma=\int_{\alpha\X_\infty}P(x,y)d\nu_\gamma=\int_{\X_\infty}(P|\alpha)(x,y)d\alpha^{-1}\nu_\gamma. \]
Since $\nu_{\alpha^{-1}\gamma\alpha}=\alpha^{-1}\nu_\gamma+\alpha^{-1}\gamma\nu_\alpha-\alpha^{-1}\nu_\alpha$, it follows that 
\[ \int_{\X_\infty}(P|\alpha)(x,y)d\alpha^{-1}\nu_\gamma=\int_{\X_\infty}(P|\alpha)(x,y)d\nu_{\alpha^{-1}\gamma\alpha}+\int_{\X_\infty}(P|\alpha)(x,y)d\alpha^{-1}(\nu_{\alpha}-\gamma\nu_\alpha). \] 
Now we show that the map $h$ on $\Gamma_0^D(Mp)$ such that $h_\gamma(P):=\int_{\X_\infty}(P|\alpha)(x,y)d\alpha^{-1}(\nu_\alpha-\gamma\nu_\alpha)$ is a coboundary. Using the fact that $\X_\infty$ is stable under the action of $\Gamma_0^D(Mp)$, an easy computation shows that 
\[ \int_{\X_\infty}(P|\alpha)(x,y)d\alpha^{-1}\gamma\nu_{\alpha}=\int_{\X_\infty}(P|\gamma\alpha)(x,y)d\alpha^{-1}\nu_{\alpha}. \]
Hence if $v(P):=-\int_{\X_\infty}(P|\alpha)(x,y)d\nu_\alpha$ then $h_\gamma(P)=(\gamma v-v)(P)$, and we are done. \end{proof}

Recall the domain of analyticity $\U_f$ introduced in Definition \ref{defi-domain}. In light of Lemma \ref{lemma-4} and normalization \eqref{normalization}, we can restate Proposition \ref{prop-pi-2} as follows.

\begin{proposition} \label{lemma-int-1}
For any choice of $\tau\in\mathcal H$ and every $k\in\U_f$ the equality 
\[ \bigg[P\longmapsto\int_{\X_{\rm aff}}P(x,y)d\tilde\mu_{f,\gamma}^\pm\bigg]=\Big[\lambda_B^\pm(k)I^\pm\bigl(f_k^{\rm JL},\tau,\gamma\bigr)\Big] \]  
holds in $H^1\bigl(\Gamma_0^D(Mp),V_{k-2}(\C)\bigr)$.
\end{proposition}

\begin{lemma} \label{lemma-int-3}
For every $k\in\U_f$, one has   
\[ \bigg[P\longmapsto\int_{\X_\infty}P(x,y)d{\tilde\mu}_{f,\gamma}^\pm\bigg]=\bigg[P\longmapsto a_p(k)^{-1}\int_{\X_{\rm aff}}P(x,y)d({\tilde\mu}_{f,\gamma}^\pm|W_p)\bigg] \]
in $H^1\bigl(\Gamma_0^D(Mp),V_{k-2}(F_k)\bigr)$.
\end{lemma}

\begin{proof} By Lemma \ref{lemma4.3}, there exists a $\Z_p$-valued measure $m$ on $\X_\infty$ such that
\[ \int_{\X_\infty}P(x,y)d{\tilde\mu}_{f,\gamma}^\pm=\int_{\X_{\rm aff}}P(x,y)d({\tilde\mu}_{f,\gamma}^\pm|U_p^{-1}W_p)+\int_{\X_{\infty}}P(x,y)d(\gamma m-m). \]
Since $\tilde\mu_f^\pm|U_p^{-1}\equiv a_p(k)^{-1}\tilde\mu_f^\pm$, there exists also an $F_k$-valued measure $m'$ on $\X_{\rm aff}$ such that 
\[ \int_{\X_{\rm aff}}P(x,y)d({\tilde\mu}_{f,\gamma}^\pm|U_p^{-1}W_p)=a_p(k)^{-1}\int_{\X_{\rm aff}}P(x,y)d(\tilde\mu_{f,\gamma}^\pm|W_p)+\int_{\X_{\rm aff}}(P|\omega_p^*)(x,y)d(\gamma m'-m') \]
where $\omega_p^*:=p\omega_p^{-1}$. Therefore, setting $v(P):=\int_{\X_\infty}Pdm+\int_{\X_{\rm aff}}(P|\omega_p^*)dm'$ we get the equality 
\[ \int_{\X_\infty}P(x,y)d{\tilde\mu}_{f,\gamma}^\pm=a_p(k)^{-1}\int_{\X_{\rm aff}}P(x,y)d(\tilde\mu_{f,\gamma}^\pm|W_p)+(\gamma v-v)(P), \] 
which completes the proof. \end{proof}

\begin{lemma} \label{lemma*} 
For any choice of $\tau\in\mathcal H$ and every $k\in\U_f$, the equality 
\[ \bigg[P\longmapsto\int_{\X_{\rm aff}}P(x,y)d\tilde\mu_{f,\gamma}^\pm|W_p\bigg]=\Big[\lambda_B^\pm(k)I^\pm\bigl(f_k^{\rm JL},\tau,\gamma\bigr)|W_p\Big] \] 
holds in $H^1\bigl(\Gamma_0^D(Mp),V_{k-2}(\C)\bigr)$. 
\end{lemma}

\begin{lemma} \label{lemma-int-2} 
For any choice of $\tau\in\mathcal H$ and every $k\in\U_f$, the equality 
\[ \bigg[P\longmapsto\int_{\X_{\infty}}P(x,y)d{\tilde\mu}_{f,\gamma}^\pm\bigg]=\bigg[P\longmapsto\frac{\lambda_B^\pm(k)p^{k-2}}{a_p(k)}I^\pm\bigl(f_k^{\rm JL},\tau,\omega_p\gamma\omega_p^{-1}\bigr)\bigl(P(x,y)|\omega_p^{-1}\bigr)\bigg] \]
holds in $H^1\bigl(\Gamma_0^D(Mp),V_{k-2}(\C)\bigr)$.
\end{lemma}

\begin{proof} Combining Lemmas \ref{lemma-int-3} and \ref{lemma*}, we get 
\[ \bigg[P\longmapsto\int_{\X_{\rm aff}}P(x,y)d({\tilde\mu}_{f,\gamma}^\pm|W_p)\bigg]=\Big[P\longmapsto\lambda_B^\pm(k)a_p(k)^{-1}I^\pm\bigl(f_k^{\rm JL},\tau,\omega_p\gamma\omega_p^{-1}\bigr)\bigl(P(x,y)|\omega_p^*\bigr)\Big]. \]
Since $\omega_p^*=p\omega_p^{-1}$, the claim follows. \end{proof}

Let $k\geq4$ be an even integer and let $f_k^\sharp$ be the $p$-stabilization of $f_k$ as defined, e.g., in \cite[p. 410]{GS}. Then $f_k^\sharp$ corresponds by Jacquet--Langands to a form $f_k^{\rm JL,\sharp}$, which can be chosen in such a way that 
\begin{equation} \label{p-stab}
f_k^{\rm JL}=f_k^{{\rm JL},\sharp}-a_p(k)^{-1}f_k^{{\rm JL},\sharp}|W_p.
\end{equation}
Combining \eqref{eq-31} and \eqref{p-stab}, for all polynomials $P$ of degree $k-2$ and any $\tau\in\mathcal H$ we get
\[ \begin{split}
   I\bigl(f_k^{\rm JL},\tau,\gamma\bigr)(P)&=\int_\tau^{\gamma(\tau)}P(z,1)f_k^{\rm JL,\sharp}(z)dz-a_p(k)^{-1}\int_\tau^{\gamma(\tau)}P(z,1)f_k^{\rm JL,\sharp}(z)|W_pdz\\
   &=I\bigl(f_k^{{\rm JL},\sharp},\tau,\gamma\bigr)(P)-a_p(k)^{-1}I\bigl(f_k^{\rm JL,\sharp}|W_p,\tau,\gamma\bigr)(P)\\
   &=I\bigl(f_k^{{\rm JL},\sharp},\tau,\gamma\bigr)(P)-p^{k-2}a_p(k)^{-1}I\bigl(f_k^{\rm JL,\sharp},\omega_p(\tau),\omega_p\gamma\omega_p^{-1}\bigr)(P|\omega_p^{-1}).
   \end{split} \] 
In particular, taking $\pm$-eigenspaces gives 
\begin{equation} \label{eq33}
I^\pm\bigl(f_k^{\rm JL},\tau,\gamma\bigr)(P)=I^\pm\bigl(f_k^{{\rm JL},\sharp},\tau,\gamma\bigr)(P)-p^{k-2}a_p(k)^{-1}I^\pm\bigr(f_k^{\rm JL,\sharp},\omega_p(\tau),\omega_p\gamma\omega_p^{-1}\bigr)(P|\omega_p^{-1}).
\end{equation} 
 
\begin{proposition} \label{prop4.9} 
Let $k\geq4$ be an even integer in $\U_f$ and let $\tau\in\mathcal H$. Then
\[ \bigg[P\longmapsto\int_\X P(x,y)d\tilde\mu_{f,\gamma}^\pm\bigg]=\bigg[\lambda_B^\pm(k)\bigg(1-\frac{p^{k-2}}{a_p(k)^2}\bigg)I^\pm\bigl(f_k^{{\rm JL},\sharp},\tau,\gamma\bigr)\bigg] \] 
in $H^1\bigl(\Gamma_0^D(Mp),V_{k-2}\C)\bigr)$.
\end{proposition}

\begin{proof} Split the integral on the left hand side  
\begin{equation} \label{eq-split}
\int_\X P(x,y)d\tilde\mu_{f,\gamma}^\pm=\int_{\X_{\rm aff}}P(x,y) d\tilde\mu_{f,\gamma}^\pm+\int_{\X_\infty} P(x,y)d\tilde\mu_{f,\gamma}^\pm.
\end{equation}
We begin by evaluating the first integral on the right hand side of \eqref{eq-split}. A straightforward computation shows that the class in $H^1\bigl(\Gamma_0^D(Mp),V_{k-2}(\C)\bigr)$ of the cocycle
\[ \gamma\longmapsto\bigg(P(x,y)\mapsto\int_\tau^{\omega_p\gamma\omega_p^{-1}(\tau)}f_k^{{\rm JL},\sharp}P(z,1)|\omega_p^{-1}dz\bigg) \]
is independent of the choice of $\tau\in\mathcal H$. Then, using Proposition \ref{lemma-int-1} and \eqref{eq33}, we see that the cocycle appearing as the first term on the right hand hand side of \eqref{eq-split} is equivalent to  
\begin{equation} \label{eq-37}
\gamma\mapsto\Big(P\mapsto\lambda_B^\pm(k)I^\pm\bigl(f_k^{{\rm JL},\sharp},\tau,\gamma\bigr)(P)-\lambda_B^\pm(k)p^{k-2}a_p(k)^{-1}I^\pm\bigl(f_{k}^{\rm JL,\sharp},\tau,\omega_p\gamma\omega_p^{-1}\bigr)(P|\omega_p^{-1})\Big).
\end{equation}
Now we compute the second term on the right of \eqref{eq-split}. Using Lemma \ref{lemma-int-2} and \eqref{eq33}, we compute 
\[ \begin{split}
   I^\pm\bigl(f_k^{\rm JL},\tau,&\,\omega_p\gamma\omega_p^{-1}\bigr)(P|\omega_p^{-1})\\
   &=I^\pm\bigl(f_k^{\rm JL,\sharp},\tau,\omega_p\gamma\omega_p^{-1}\bigr)(P|\omega_p^{-1})-p^{k-2}a_p(k)^{-1}I^\pm\bigl(f_k^{\rm JL,\sharp},\omega_p(\tau),\omega_p^2\gamma\omega_p^{-2}\bigr)(P|\omega_p^{-2}).
   \end{split} \] 
To study the second term on the right hand side of the above equation, we apply \eqref{eq-31}, the equality $f_k^{{\rm JL},\sharp}|W_p^2=p^{k-2}f_k^{{\rm JL},\sharp}$ and the fact that the class of $\gamma\mapsto I^\pm\bigl(f_k^{\rm JL,\sharp},\tau,\gamma\bigr)$ does not depend on the choice of $\tau\in\mathcal H$. Finally, we obtain that $P\mapsto\int_{\X_\infty}P(x,y)d\tilde\mu_{f,\gamma}^\pm$ is equivalent to 
\begin{equation} \label{eq39}
P\longmapsto\lambda_B^\pm(k)a_p(k)^{-1}p^{k-2}\Big(I^\pm\bigl(f_k^{\rm JL,\sharp},\tau,\omega_p\gamma\omega_p^{-1}\bigr)(P|\omega_p^{-1})-a_p(k)^{-1}I^\pm\bigl(f_k^{\rm JL,\sharp},\tau,\gamma\bigr)(P)\Big).
\end{equation}
Combining \eqref{eq-37} and \eqref{eq39} yields the result. \end{proof}

\begin{corollary} \label{coro4.9} 
Let $k\geq4$ be an even integer in $\U_f$ and let $\tau\in\mathcal H$. Then
\[ \bigg[P\longmapsto\int_\X P(x,y)d\tilde\mu_{f,\gamma}^\pm\bigg]=\bigg[\lambda_B^\pm(k)\bigg(1-\frac{p^{k-2}}{a_p(k)^2}\bigg)I^\pm\bigl(f_k^{{\rm JL},\sharp},\tau ,\gamma\bigr)\bigg] \] 
in $H^1\bigl(\Gamma_0^D(M),V_{k-2}(\C)\bigr)$.
\end{corollary}

\begin{proof} If this were not true then the left hand side and the right hand side would represent different classes in $H^1\bigl(\Gamma_0^D(M),V_{k-2}(\C)\bigr)$ 
that, by Proposition \ref{prop4.9}, when restricted to $\Gamma_0^D(Mp)$ share the same Hecke eigenvalues at all primes except, possibly, the prime $p$. By multiplicity one, this cannot occur because $f_k^\sharp$ is a newform. \end{proof}

\begin{proposition} \label{prop-split} 
For any $\tau\in\mathcal H$, the equality 
\[ \bigg[\gamma\longmapsto\int_{\Z_p^\times\times\Z_p^\times}d\tilde\mu_{f,\gamma}^\pm\bigg]=\bigg[\gamma\longmapsto2\lambda_B^\pm(2)\cdot\int_\tau^{\gamma(\tau)}f^{\rm JL}dz\bigg] \]
holds in $H^1\bigl(\Gamma_0^D(Mp),\C\bigr)$. 
\end{proposition}

\begin{proof} Combining \eqref{eq-31}, Proposition \ref{lemma-int-1} and Lemma \ref{lemma-int-2}, it follows that 
\[ \int_{\Z_p^\times\times\Z_p^\times}d\tilde\mu_{f,\gamma}^\pm=\lambda_B^\pm(2)\bigg(\int_\tau^{\gamma(\tau)}f^{\rm JL}dz-\frac{1}{a_p}\int_\tau^{\gamma(\tau)}(f^{\rm JL}|W_p)dz\bigg) \]
for all $\gamma\in\Gamma_0^D(Mp)$. Now $f^{\rm JL}|W_p=-a_p f^{\rm JL}$, and the result follows. \end{proof}

Define $w_p:=-a_p$, so that $f_k^{{\rm JL},\sharp}|W_p=w_pp^{(k-2)/2}f_k^{{\rm JL},\sharp}$ (it is known that the value $w_p$ is common to all the forms $f_k$ in the Hida family: see, e.g., \cite[\S 3.4.4]{NP}). 

\begin{proposition} \label{prop-split-2}
Let $k\geq 4$ be an even integer in $\U_f$ and let $\tau\in\mathcal H$. Then
\[ \bigg[P\longmapsto\int_{\Z_p^\times\times\Z_p^\times}P(x,y)d\tilde\mu_{f,\gamma}^\pm\bigg]=\bigg[\lambda_B^\pm(k)\bigg(1-\frac{w_pp^{(k-2)/2}}{a_p(k)}\bigg)^2I^\pm\bigl(f_k^{{\rm JL},\sharp},\tau,\gamma\bigr)\bigg] \]
in $H^1\bigl(\Gamma_0^D(Mp),V_{k-2}(\C)\bigr)$.
\end{proposition}

\begin{proof} Combining Lemma \ref{lemma-4} with the arguments in the proof of Proposition \ref{prop4.9}, it follows that the cocyle on $\Gamma_0^D(Mp)$ sending $\gamma$ to $\bigl(P\mapsto \int_{\Z_p^\times\times\Z_p^\times}Pd\tilde\mu_{f,\gamma}^\pm\bigr)$ is equivalent to the difference between the cocycles sending $\gamma$ to the maps which take $P$ to 
\[ \lambda_B^\pm(k)I^\pm\bigl(f_k^{{\rm JL},\sharp},\tau,\gamma\bigr)(P)-\frac{\lambda_B^\pm(k)p^{k-2}}{a_p(k)}I^\pm\bigl(f_k^{\rm JL,\sharp},\tau,\omega_p\gamma\omega_p^{-1}\bigr)(P|\omega_p^{-1}) \] 
and 
\[ \frac{\lambda_B^\pm(k)p^{k-2}}{a_p(k)}I^\pm\bigl(f_k^{\rm JL,\sharp},\tau,\omega_p\alpha^{-1}\gamma\alpha\omega_p^{-1}\bigr)(P|\alpha\omega_p^{-1})-\frac{\lambda_B^\pm(k)p^{k-2}}{a_p(k)^2}I^\pm\bigl(f_k^{\rm JL,\sharp},\tau,\alpha^{-1}\gamma\alpha\bigr)(P|\alpha), \] 
with $\alpha$ as in the statement of Lemma \ref{lemma-4}. Using \eqref{eq-31} and the fact that $f_k^{\rm JL,\sharp}$ and $f_k^{\rm JL,\sharp}|W_p$ have level $\Gamma_0^D(M)$, it can be checked that the sums above are equal to 
\[ \lambda_B^\pm(k)I^\pm\bigl(f_k^{{\rm JL},\sharp},\tau,\gamma\bigr)(P)-\frac{\lambda_B^\pm(k)}{a_p(k)}I^\pm\bigl(f_k^{\rm JL,\sharp}|W_p,\tau,\gamma\bigr)(P) \]
and
\[ \frac{\lambda_B^\pm(k)}{a_p(k)}I^\pm\bigl(f_k^{\rm JL,\sharp}|W_p,\tau,\gamma\bigr)(P)-\frac{\lambda_B^\pm(k)p^{k-2}}{a_p(k)^2}I^\pm\bigl(f_k^{\rm JL,\sharp},\tau,\gamma\bigr)(P),\]
respectively. Taking the difference of these two expressions yields the result. \end{proof}

\begin{corollary} \label{coro-split}
Let $k\geq4$ be an even integer in $\U_f$ and let $\tau\in\mathcal H$. Then
\[ \bigg[P\longmapsto\int_{\Z_p^\times\times\Z_p^\times}P(x,y)d\tilde\mu_{f,\gamma}^\pm\bigg]=\bigg[\lambda_B^\pm(k)\bigg(1-\frac{w_pp^{(k-2)/2}}{a_p(k)}\bigg)^{\!\!2}I^\pm\bigl(f_k^{{\rm JL},\sharp},\tau,\gamma\bigr)\bigg] \]
in $H^1\bigl(\Gamma_0^D(M),V_{k-2}(\C)\bigr)$.
\end{corollary}

\begin{proof} In light of Proposition \ref{prop-split-2}, proceed as in the proof of Corollary \ref{coro4.9}. \end{proof}

\subsection{Popa's work} \label{popa-sec}

Let $\mathbb A$ be the adele ring of $\Q$ and set $B_\mathbb A:=B\otimes_\Q\mathbb A$ for the adelization of our quaternion algebra $B$. Let $g$ be a weight $k$ normalized newform of level $\Gamma_0({S} D)$ where ${S} $ is a square-free 
positive integer prime to $D$. Let $R_0^D({S} )$ be an Eichler order of $B$ of level ${S} $ and let $\Gamma_0^D({S} )$ denote the group of its elements of norm $1$. Let $g^{\rm JL}$ denote a modular form on $\Gamma_0^D({S} )$ corresponding to $g$ via the Jacquet--Langlands correspondence; notice that $g^{\rm JL}$ is well defined up to multiplication by non-zero constants. Let $\varphi_g$ and $\varphi_g^{\rm JL}$ denote the automorphic forms on $\GL_2(\mathbb A)$ and $B_\mathbb A^\times$ associated with $g$ and $g^{\rm JL}$, respectively, whose definitions can be found, e.g., in \cite[\S 3]{Gel} and \cite[\S 10]{Gel}, respectively. We just recall, because it will be used in the following computations, that $\varphi_g^{\rm JL}$ is defined in \emph{loc. cit.} using the decomposition 
\begin{equation} \label{quaternion-adelic-decomposition}
B_\mathbb A^\times=B^\times(B_\infty^\times)^0K_0^D({S} )
\end{equation}
where $K_0^D({S} ):=\big(R_0^D({S} )\otimes_\Z\hat\Z\big)^\times$ ($\hat\Z$ denotes, as usual, the profinite completion of $\Z$) and $(B_\infty^\times)^0$ is the connected component of the identity in $B_\infty^\times:=(B\otimes_\Q\R)^\times$ (cf. \cite[Lemma 10.3]{Gel} for details). 

Fix a real quadratic field $L=\Q(\sqrt{d_L})$ such that all the primes dividing ${S} $ (respectively, $D$) are split (respectively, inert) in $L$. For simplicity, we assume that $L\neq F$ (actually, the case $L=F$ is easier but requires different conventions: note, for example, that one may diagonally embed $F$ into $\M_2(F)$ and that the quadratic form associated with this embedding has $c=b=0$, in the notations of \S\ref{sec4.1-bis}). Set $G_L^+:=\Gal(H_L^+/L)$ where $H_L^+$ is the narrow Hilbert class field of $L$. Fix an unramified character $\chi:G_L^+\rightarrow\C^\times$. By class field theory, $\chi$ can be identified with a finite-order Hecke character $\chi:\mathbb A^\times L^\times\backslash\mathbb A_L^\times\rightarrow\C^\times$, where $\mathbb A_L=L\otimes_\Q\mathbb A$ is the adele ring of $L$. Denote by $\pi_\chi$ the automorphic representation of $\GL_2(\mathbb A)$ attached to $\chi$ via the Jacquet--Langlands correspondence. 

Fix an embedding $L\hookrightarrow\R$ and choose the generator $\varepsilon_L$ of the group of norm $1$ elements in $\cO_L$ so that $\varepsilon_L>1$ with respect to this embedding. In the special case where $L=K$ fix the embedding $K\hookrightarrow\R$ as in \S \ref{sec2.1}, so that $\varepsilon_L$ is the element denoted by $\varepsilon_1$ in \S\ref{sec2.5}. Write $dx$ for the Haar measure on $L^\times\mathbb A^\times\backslash\mathbb A^\times_L$ defined in \cite[p. 840]{Po}, which is normalized in such a way that the total mass is equal to $h^+_L\ln\varepsilon_L$, where $h^+_L$ is the cardinality of $\Pic^+(\mathcal O_L)$. As above, write $\Emb\bigl(\mathcal O_L,R_0^D({S} )\bigr)$ for the set of optimal embeddings of $\mathcal O_L$ into $R_0^D({S} )$. Fix also an embedding $\psi_0\in\Emb\bigl(\mathcal O_L,R_0^D({S} )\bigr)$ (if $L=K$ choose $\psi_0$ as in \S \ref{Galois-subsec}) and write $\psi_{0,\mathbb A}$ for its adelization. Also, write $\hat\psi_0$ and $\psi_{0,\infty}$ for the finite and the infinite part of $\psi_{0,\mathbb A}$, respectively. At the archimedean places of $L$ we may also consider, as in \cite[(5.3.1)]{Po}, the diagonal embedding $\psi_{\rm diag}:L\otimes_\Q\R\hookrightarrow\M_2(\R)$ obtained by extending by $\R$-linearity the embedding $L\hookrightarrow\M_2(\R)$ given by $x\mapsto\smallmat x00{\bar x}$ where $x\mapsto\bar x$ is the involution of $L$. Let $i_F\big(\psi_0(\sqrt{d_L})\big)=\smallmat abc{-a}$ and define $\gamma_\infty:=\smallmat{a+\sqrt{d_L}}{a-\sqrt{d_L}}{c}{c}\in\GL_2(\R)$ with $c>0$ (note that, since $L\neq F$, we always have $c\neq0$). Then $\psi_{0,\infty}=\gamma_\infty\psi_{\rm diag}\gamma_\infty^{-1}$. Define 
\[ l(\varphi_g^{\rm JL},\chi):=\int_{\mathbb A^\times L^\times\backslash\mathbb A_L^\times}\varphi_g^{\rm JL}\bigl({\psi_{0,\mathbb A}}(x)\gamma_\infty\bigr)
\chi^{-1}(x)dx. \]
By \cite[Theorem 5.4.1]{Po}, there is the formula 
\[ L(1/2,\pi_g\times\pi_\chi)=\frac{4}{\sqrt{d_L}}\cdot\prod_{\ell|D}\frac{\ell+1}{\ell-1}\cdot\frac{\|\varphi_g\|^2}{\|\varphi_g^{\rm JL}\|^2}\cdot|l(\varphi_g^{\rm JL},\chi)|^2. \] 
Although $l(\varphi_g^{\rm JL},\chi)$ depends on the choice of $g^{\rm JL}$ and $\varphi_g^{\rm JL}$, the right hand side of the above expression does not depend on it.

Now we compute the term $l(\varphi_g^{\rm JL},\chi)$ more explicitly. Set $\tau_0:=i_F(\gamma_\infty)(i)$. Following \cite[\S 6.1]{Po}, one shows that 
\[ l(\varphi_g^{\rm JL},\chi)=\frac{1}{(2i)^{k/2}\sqrt{d_L}^{(k-2)/2}}\sum_{[\mathfrak a]\in\Pic^+(\mathcal O_L)}\chi^{-1}([\fr a])M([\fr a]) \]
with 
\begin{equation} \label{M-formula-eq}
M([\fr a]):=\det(\alpha_{[\fr a]})^{\frac{2-k}{2}}\int_{\tau_{[\fr a]}}^{(\alpha_{[\fr a]}^{-1}\psi_0\alpha_{[\fr a]})(\varepsilon_L)\tau_{[\fr a]}}f(z)\bigl(Q_{\psi_0}^{(k)}|\alpha_{[\fr a]}\bigr)(z,1)dz.
\end{equation}
The term $\alpha_{[\fr a]}\in(B_\infty^\times)^0$ is defined, using decomposition \eqref{quaternion-adelic-decomposition}, by the equation $\hat\psi_0(\hat{\fr a})=b\alpha_{[\fr a]}^{-1}\boldsymbol k$ with $b\in B^\times$, $\alpha_{[\fr a]}^{-1}\in (B_\infty^\times)^0$, $\boldsymbol k\in {{K}}_0^D({S} )$ and $\tau_{[\fr a]}:=\alpha_{[\fr a]}^{-1}(\tau_0)$; here $\hat{\fr a}$ is any finite idele that corresponds to an ideal $\fr a$ representing the class $[\fr a]$. Although $\alpha_{[\fr a]}$ is not unique, the results will not depend on the actual choice. To simplify notations, we also identify $\alpha_{[\fr a]}$ with its image via $i_\infty$ and view $\alpha_{[\fr a]}$ as an element of $\GL_2^+(\R)$.

Let $h_L$ be the class number of $L$, so that $h_L^+/h_L=1$ or $2$, let $\Pic(\cl O_L)$ be the ideal class group of $L$ and write $H_L$ for the Hilbert class field of $L$ (hence $\Pic(\cO_L)\simeq\Gal(H_L/L)$ via the reciprocity map of class field theory). Let $\mathfrak D_L$ denote the class of $(\sqrt{d_L})$ in $\Pic^+(\cO_L)$. If for every $\sigma\in \Gal(H_L/L)$ we write $[\fr a_\sigma]=\text{rec}^{-1}(\sigma)$ for an ideal $\fr a_\sigma$ of $\cl O_L$ then 
\[ \Pic(\cl O_L)=\bigl\{[\fr a_\sigma]\bigr\}_{\sigma\in\Gal(H_L/L)} \]
and
\[ \Pic^+(\cl O_L)=\begin{cases}\bigl\{[\fr a_\sigma]\bigr\}_{\sigma\in\Gal(H_L/L)}&\text{if $h_L^+=h_L$},\\[2mm]\bigl\{[\fr a_\sigma]\bigr\}_{\sigma\in\Gal(H_L/L)}\cup\bigl\{\fr D_L\cdot[\fr a_\sigma]\bigr\}_{\sigma\in\Gal(H_L/L)}&\text{otherwise}.\end{cases} \]
Here, with a slight abuse of notation, we have adopted the symbol $[\star]$ to denote ideal classes both in $\Pic(\cl O_L)$ and in $\Pic^+(\cl O_L)$.

Write $\sigma_L\in G_L^+$ for the image of $\fr D_L$ under the reciprocity map. For every $\sigma\in G_L^+$ let $[\mathfrak a_\sigma]$ be the corresponding element in $\Pic^+(\mathcal O_L)$ via the Artin map, represented by an ideal $\fr a_\sigma\subset\cO_L$. We extend in the obvious way the notation and results of \S \ref{Galois-emb-subsec} about oriented optimal embeddings to the more general case needed here. In particular, for an oriented optimal embedding $\psi\in\mathcal E\big(\cO_L,R_0^D({S} )\big)$ define 
\begin{equation} \label{eq70}
[\fr a_\sigma]\star[\psi]:=\bigl[a_\sigma\psi a_\sigma^{-1}\bigr]\in\mathcal E\bigl(\cO_L,R_0^D(S)\bigr)/\Gamma_0^D(S)
\end{equation}
where $a_\sigma\in R_0^D(S)$ has positive reduced norm and satisfies $R_0^D(S)\psi(\mathfrak a_\sigma)=R_0^D(S)a_\sigma$. Fixing $\psi$, the rule $\sigma\mapsto [\psi_\sigma]$ with $\psi_\sigma:= a_\sigma\psi a_\sigma^{-1}$ is the inverse of a bijection $G$ between $\mathcal E\big(\cO_L,R_0^D(S)\big)/\Gamma_0^D(S)$ and $G_L^+$. The map $G$ satisfies the equation $G([\psi^*])=\sigma_L\cdot G([\psi])$ for all $\psi\in\E\bigl(\cO_L,R_0^D(S)\bigr)$, where $\psi^*:=\omega_\infty\psi\omega_\infty^{-1}$. The family $\{\psi_\sigma\}_{\sigma\in G^+_L}$ is a set of representatives of the $\Gamma_0^D(S)$-conjugacy classes of oriented optimal embeddings 
of $\cl O_L$ into $R_0^D(S)$. Finally, for every $\sigma\in G_L^+$ set $\gamma_\sigma:=\psi_\sigma(\varepsilon_L)\in\Gamma_0^D(S)$.

\begin{remark} \label{remark4.16}
We explicitly observe a fact that will be used in the proof of Theorem \ref{thm-factor}. Suppose that, in the above discussion, $S=Mp$. Since the quadratic order $\mathcal O_L$ is maximal, an oriented optimal embedding of $\mathcal O_L$ into $R_0^D(Mp)$ can be naturally viewed as an element of $\mathcal E\bigl(\mathcal O_L,R_0^D(M)\bigr)$ too. Moreover, the element $a_\sigma$ described above also satisfies $R_0^D(M)\psi(\mathfrak a_\sigma)=R_0^D(M)a_\sigma$. This shows that we can (and do) fix via \eqref{eq70} a family $\{\psi_\sigma\}_{\sigma\in G^+_L}$ which is a set of representatives of the $\Gamma_0^D(S)$-conjugacy classes of oriented optimal embeddings of $\cl O_L$ into $R_0^D(S)$ for $S$ equal to either $M$ or $Mp$.
\end{remark}

If the subscript $\infty$ indicates that an element of $B^\times$ of positive norm is to be viewed as belonging to $(B^\times_\infty)^0$, it follows that $\hat\psi_0(\hat{\fr a}_\sigma)=a_\sigma{(a^{-1}_\sigma)}_\infty\boldsymbol k$, with respect to decomposition \eqref{quaternion-adelic-decomposition}, for a suitable $\boldsymbol k\in K_0^D(S)$. Hence for every $\sigma\in G_L^+$ we can take $\alpha_{[\fr a_\sigma]}=a_\sigma$, in the above notation. Moreover, if $h_L^+\not=h_L$ then for every $\sigma\in G^+_L$ we can also take $\alpha_{\fr D_L\cdot[\fr a_\sigma]}$ to be equal to $\omega_\infty a_\sigma\psi_0(\sqrt{d_L})$.  

In light of the expression of the integrals $M([\fr a_\sigma])$ given in \eqref{M-formula-eq}, we summarize the above calculations in the following 
\begin{theorem}[Popa] \label{teo-popa} 
Let $\chi$ be a complex-valued character of $G_L^+$. Then 
\[ L(1/2,\pi_g\times\pi_\chi)=\frac{4}{(2i)^k\sqrt{d_L}^{k-1}}\cdot\frac{\|\varphi_g\|^2}{\|\varphi_g^{\rm JL}\|^2}\cdot\prod_{\ell|D}\frac{\ell+1}{\ell-1}\cdot\Bigg|\sum_{\sigma\in G_L^+}\chi^{-1}(\sigma)\int_{\tau_0}^{\gamma_\sigma(\tau_0)}g^{\rm JL}(z)Q_{\psi_\sigma}^{(k)}(z,1)dz\Bigg|^2. \]
\end{theorem}

With notation as in \eqref{pm-I-short-eq}, the next task is to rewrite more explicitly the integral
\[ \Theta_\psi:=\int_{\tau_0}^{\gamma_\psi(\tau_0)}g^{\rm JL}(z)Q_\psi^{(k)}(z,1)dz=I\bigl(g^{\rm JL},\psi\bigr) \] 
for any of the $\psi=\psi_\sigma$ and $\gamma_\psi=\gamma_\sigma$. Define 
\[ (g^{\rm JL}|\omega_\infty)(z):=(Cz+D)^{-k}\cdot\overline{g^{\rm JL}(\omega_\infty\bar z)}, \] 
where $i_\infty(\omega_\infty)=\smallmat **CD$. 

\begin{lemma} \label{rescale}
Up to rescaling $g^{\rm JL}$, we may assume that $g^{\rm JL}|\omega_\infty=g^{\rm JL}$. 
\end{lemma}

\begin{proof} A simple computation shows that $g^{\rm JL}|\omega_\infty$ has the same weight and level as $g^{\rm JL}$. Since the Hecke eigenvalues of $g$ are real because $g$ is a newform with trivial character, one sees that the same is true of $g^{\rm JL}$. This can be used to show that $g^{\rm JL}|\omega_\infty$ has the same eigenvalues as $g^{\rm JL}$ and so, by multiplicity one, there exists a complex number $\lambda$ such that $g^{\rm JL}|\omega_\infty=\lambda g^{\rm JL}$. Since $\omega_\infty^2\in\Gamma_0^D(M)$, we have that $\bigl(g^{\rm JL}|\omega_\infty\bigr)|\omega_\infty=g^{\rm JL}$. (Thus, in particular, $|\lambda|=1$, but we will not need this here). If $\lambda=1$ then $g^{\rm JL}$ verifies the claimed condition without rescaling, while if $\lambda=-1$ then $ig^{\rm JL}$ does the job. In the other cases, the modular form $g^{\rm JL}+g^{\rm JL}|\omega_\infty=(1+\lambda)g^{\rm JL}$ has the required property. \end{proof}

From now on assume that $g^{\rm JL}$ verifies the condition in Lemma \ref{rescale}. Using this, a simple computation shows that
\[ \overline{\Theta_\psi}=\int_{\omega_\infty\bar\tau_0}^{\gamma_{\psi^*}(\omega_\infty\bar\tau_0)}g^{\rm JL}(z)Q^{(k)}_{\psi^*}(z,1)dz=I\bigl(g^{\rm JL},\psi^*\bigr)=\Theta_{\psi^*} \] 
(recall that $\psi^*:=\omega_\infty\psi\omega_\infty^{-1}$). By the above discussion, the embeddings $\psi_\sigma^*$ and $\psi_{\sigma_L\sigma}$ are $\Gamma_0^D(S)$-conjugate. Since $\Theta_\psi$ depends only on the $\Gamma_0^D(S)$-conjugacy class of $\psi$, we obtain 
\begin{equation} \label{eq-Theta-2}
\overline{\Theta_{\psi_\sigma}}=\Theta_{\psi_{\sigma_L\sigma}}.
\end{equation}
Let $\chi$ be a \emph{genus character} of $L$, i.e., an unramified quadratic character of $G_L^+$. Thanks to \eqref{eq-Theta-2}, we have 
\[ \overline{\sum_{\sigma\in G_L^+}\chi(\sigma)\Theta_{\psi_\sigma}}=\sum_{\sigma\in G_L^+}\chi(\sigma)\overline{\Theta_{\psi_\sigma}}=\sum_{\sigma\in G_L^+}\chi(\sigma)\Theta_{\psi_{\sigma_L\sigma}}=\chi(\sigma_L)\cdot\Bigg(\sum_{\sigma\in G_L^+}\chi(\sigma)\Theta_{\psi_\sigma}\Bigg). \]

\begin{definition}
The \emph{sign} of $\chi$ is the sign of $\chi(\sigma_L)$.
\end{definition}
 
A straightforward calculation shows that
\[ \overline{\int_{\tau_0}^{\gamma(\tau_0)}g(z)P(z,1)dz}=-\int_{\omega_\infty\bar\tau_0}^{\omega_\infty\gamma\omega_\infty^{-1}(\omega_\infty\bar\tau_0)}g(z)(P|\omega_\infty^*)(z,1)dz. \] 
In other words, the action of complex conjugation on $I(g,\tau,\gamma)(P)$ coincides, up to a change of sign, with the action of $W_\infty$ (for a more detailed discussion of this relation in the case of elliptic modular forms, see \cite[p. 588]{hida}). Hence if $\chi$ has sign $-\epsilon\in\{\pm\}$ then 
\[ \sum_{\sigma\in G_L^+}\chi(\sigma)\Theta_{\psi_\sigma}=\sum_{\sigma\in G_L^+}\chi(\sigma)I^\epsilon\bigl(g^{\rm JL},\psi_\sigma\bigr). \] 
Summing up, Theorem \ref{teo-popa} admits the following reformulation for genus characters.

\begin{corollary} \label{coro-popa} 
Let $\chi$ be a genus character of $L$ of sign $-\epsilon$. Then 
\[ L(1/2,\pi_g\times\pi_\chi)=\frac{\epsilon4}{(2i)^k\sqrt{d_L}^{k-1}}\cdot\frac{\|\varphi_g\|^2}{\|\varphi_g^{\rm JL}\|^2}\cdot\prod_{\ell|D}\frac{\ell+1}{\ell-1}\cdot\Bigg(\sum_{\sigma\in G_L^+}\chi(\sigma)I^\epsilon\bigl(g^{\rm JL},\psi_\sigma\bigr)\Bigg)^{\!\!2}. \]
\end{corollary}

\subsection{$p$-adic $L$-functions over real quadratic fields} 

Before studying $p$-adic $L$-functions, let us observe a simple fact. Let $L$ be a real quadratic field and let $\psi:L\hookrightarrow B$ be an optimal embedding of $\mathcal O_L$ into $R_0^D(N)$, where $N$ denotes either $M$ or $Mp$. Fix a subset $\U$ of $\X$ such that $\gamma_\psi\U=\U$. Let $k\geq2$ be an \emph{even} integer, fix $\boldsymbol\nu\in\W$ and choose a representative $\nu$ of $\boldsymbol\nu$. 
\begin{lemma} \label{well-def} 
Let the notation be as above. 
\begin{enumerate}
\item The value $\int_\U Q^{(k)}_\psi(x,y)d\nu_{\gamma_\psi}$ does not depend on the choice of a representative $\nu$ of $\boldsymbol\nu$.  
\item If $\psi'=\alpha\psi\alpha^{-1}$ for some $\alpha\in\Gamma_0^D(N)$ then $\int_\U Q^{(k)}_\psi(x,y)d\nu_{\gamma_\psi}=\int_{\alpha\U}Q^{(k)}_{\psi'}(x,y)d\nu_{\gamma_\psi}$.
\item If $\gamma\U=\U$ for all $\gamma\in\Gamma_0^D(N)$ then $\int_\U Q^{(k)}_\psi(x,y)d\nu_{\gamma_\psi}$ depends only on the $\Gamma_0^D(N)$-conjugacy class of $\psi$.
\end{enumerate}
\end{lemma} 

\begin{proof} The proof follows from a direct calculation using the equality 
\begin{equation} \label{extra}
Q_\psi\bigl(\gamma_\psi(x,y)\bigr)=Q_\psi(x,y),
\end{equation} 
which is a consequence of \eqref{eq-21}. Details are left to the reader. \end{proof}

Recall now the fixed real quadratic field $K$ where $p$ is inert. Recall also the domain of analiticity $\U_f$ of $\boldsymbol{\tilde\mu}_f^\pm$ introduced in Definition \ref{defi-domain}. Let $\langle x\rangle$ denote the principal unit attached to $x\in\Q_p^\times$, defined as the unique element of $1+p\Z_p$ such that $x=p^{\ord_p(x)}\zeta_x\langle x\rangle$ where $\zeta_x$ is a $(p-1)$-st root of unity. For $k\in\U_f$ and $(x,y)\in\mathbb X$ define 
\[ Q_\psi^{(k)}(x,y):=Q_\psi(x,y)^\frac{k-2}{2}=\exp_p\bigg(\frac{k-2}{2}\log\bigl(\big\langle Q_\psi(x,y)\big\rangle\bigr)\!\bigg) \] 
where $\exp_p$ is the $p$-adic exponential and $\log$ is \emph{any} branch of the $p$-adic logarithm. Since $Q_\psi(x,y)=\langle Q_\psi(x,y)\rangle$ for $(x,y)\in\mathbb X$, this definition does not depend on the choice of $\log$.

\begin{definition} \label{defi-partial} 
Let $\psi$ be an optimal embedding of $\mathcal O_K$ into $R_0^D(M)$. The \emph{partial square root $p$-adic $L$-function} attached to $f_\infty$, $K$ and $\psi$ is 
\[ \mathcal L_p^\pm(f_\infty,\psi,k):=\int_\X Q_\psi^{(k)}(x,y)d{\tilde\mu}^\pm_{f,{\gamma_\psi}}(x,y) \] 
where $k$ varies in $\U_f$.
\end{definition}

\begin{remark} \label{rem4.24} 
Since $Q_\psi(x,y)=\langle Q_\psi(x,y)\rangle$ on $\mathbb X$, we have   
\[ \mathcal L_p^\pm(f_\infty,\psi,k)=\int_\X\langle Q_\psi(x,y)\rangle^{\frac{k-2}{2}}d{\tilde\mu}^\pm_{f,{\gamma_\psi}}(x,y). \]
Moreover, if $k\geq2$ is an integer in $\U_f$ then $Q_\psi^{(k)}(x,y)$ is equal to the $(k-2)/2$-fold self product of $Q_\psi(x,y)$. This justifies the restriction of $\U_f$ to the residue class of $2$ modulo $p-1$ in $\mathcal X$.
\end{remark}

\begin{proposition} \label{prop2.24}
$\mathcal L_p^\pm(f_\infty,\psi,{k} )$ is independent of the choice of a representative $\tilde\mu_f^\pm$ of $\boldsymbol{\tilde{\mu}}_f^\pm$ and only depends on the $\Gamma_0^D(M)$-conjugacy class of $\psi$. 
\end{proposition} 

\begin{proof} Lemma \ref{well-def} shows that the restriction of $\mathcal L_p^\pm(f_\infty,\psi,k)$ to the subset of $\U_f$ consisting of points with trivial character has the required properties. The result then follows from a density argument. \end{proof}

\begin{proposition} \label{prop4.15} 
For every even integer $k\geq4$ in $\U_f$ one has   
\[ \mathcal L^\pm(f_\infty,\psi,k)=\lambda_B^\pm(k)\bigg(1-\frac{p^{k-2}}{a_p(k)^2}\bigg)I^\pm\bigl(f_k^{{\rm JL},\sharp},\psi\bigr). \]
\end{proposition}
\begin{proof} Corollary \ref{coro4.9} ensures the existence of $v\in V_{k-2}(\C)$ such that 
\[ \int_\X Q_\psi^{(k)}(x,y)d{\tilde\mu}^\pm_{f,\gamma_\psi}(x,y)=\lambda_B^\pm(k)\bigg(1-\frac{p^{k-2}}{a_p(k)^2}\bigg)I^\pm\bigl(f_k^{{\rm JL},\sharp},\psi\bigr)+(\gamma v-v)\Big(Q_\psi^{(k)}(x,y)\Big) \] 
for all $\gamma\in\Gamma_0^D(Mp)$. We conclude by evaluating this expression at $\gamma=\gamma_\psi$ and using \eqref{extra}. \end{proof} 

\begin{definition} \label{defi-L_p}
Let $\chi$ be a quadratic character of $G_K^+$ of sign $-\epsilon$ and let $k$ vary in $\U_f$. 
\begin{enumerate}
\item The \emph{square root $p$-adic $L$-function} attached to $f_\infty$ and $\chi$ is
\[ \mathcal L_p(f_\infty/K,\chi,k):=\sum_{\sigma\in G_K^+}\chi(\sigma)\mathcal L_p^\epsilon(f_\infty,\psi_\sigma,k). \] 
\item The \emph{$p$-adic $L$-function} attached to $f_\infty$ and $\chi$ is 
\[ L_p(f_\infty/K,\chi,k):=\mathcal L_p(f_\infty/K,\chi,k)^2. \]
\end{enumerate}
\end{definition}

\begin{theorem} \label{teo4.17} 
Let $\chi$ be a genus character of $K$ of sign $-\epsilon$. Then for all even integers $k\geq4$ in $\U_f$ we have  
\[ L_p(f_\infty/K,\chi,k)=\lambda_B^\epsilon(k)^2\bigg(1-\frac{p^{k-2}}{a_p(k)^2}\bigg)^{\!\!2}\prod_{\ell|D}\frac{\ell-1}{\ell+1}\cdot\frac{(2i)^{k}\sqrt{d_K}^{k-1}}{\epsilon4}\cdot\frac{\|\varphi_{f_k}^{\rm JL,\sharp}\|^2}{\|\varphi_{f_k^\sharp}\|^2}\cdot L_K\bigl(1/2,\pi_{f_k^\sharp}\times\pi_\chi\bigr). \]
\end{theorem}
\begin{proof} If $k\geq2$ then $f_k^\sharp$ is a newform on $\Gamma_0^D(M)$, so the result follows immediately from Corollary \ref{coro-popa} and Proposition \ref{prop4.15}. \end{proof}

\begin{remark}
The above result holds also for $k=2$ if we set $f_2^\sharp:=0$. Indeed, in this case 
\begin{equation} \label{vanishing-eq}
\mathcal L_p(f_\infty/K,\chi,2)=\int_\X d\tilde\mu_{f,\gamma_\psi}=\int_{\PP^1(\Q_p)}d\pi_*(\mu_{f,\gamma_\psi})=0,
\end{equation}
and the equality in the statement is trivially verified. 
\end{remark}

\subsection{Atkin--Lehner involutions and genus characters} \label{atkin-genus-subsec}

A genus character $\chi$ of $K$ cuts out the genus field $H_\chi$ of $K$ given by
\[ H_\chi=\Q\bigl(\sqrt{d_1},\sqrt{d_2}\bigr) \]
where $d_K=d_1d_2$. The extension $H_\chi/\Q$ is biquadratic unless $\chi$ is trivial (in which case, of course, $H_\chi=K$). Let $\chi_1$, $\chi_2$ and $\epsilon_K$ be the characters associated with the quadratic fields $\Q(\sqrt{d_1})$, $\Q(\sqrt{d_2})$ and $K$, respectively, so that $\epsilon_K=\chi_1\cdot\chi_2$. In fact, the genus characters of $K$ are in bijection with the factorizations of $d$ into a product of relatively prime discriminants $d_1$ and $d_2$ or, equivalently, with the unordered pairs $(\chi_1,\chi_2)$ of primitive quadratic Dirichlet characters of coprime conductors satisfying $\epsilon_K=\chi_1\cdot\chi_2$ (the trivial character corresponds to the factorization $d=1\cdot d$). For more details see, e.g., \cite[Ch. 14, \S G]{Cohn}.

Let now $\psi$ be an oriented optimal embedding of $\cO_K$ into $R_0^D(M)$ and recall the unique element $\sigma_\psi\in G^+_K$ introduced in \eqref{W} such that
\begin{equation} \label{sigma-psi-eq}
\psi\circ\tau\sim\omega_{MD}(\sigma_\psi\star\psi)\omega_{MD}^{-1}. 
\end{equation}
To ease the writing, from now until the end of this subsection set $G:=G_K^+$. Denote by $G^2$ the subgroup of $G$ consisting of the squares of the elements of $G$. 

\begin{lemma} \label{square-independent-lemma}
The class of $\sigma_\psi$ in $G/G^2$ does not depend on $\psi$.
\end{lemma}

In light of this result, write $\bar\sigma$ for the class of any $\sigma_\psi$ in $G/G^2$.

\begin{proof} To begin with, recall that $\sigma_\psi$ is characterized as the unique element of $G$ such that
\[ \bigl[\omega_{MD}^{-1}(\psi\circ\tau)\omega_{MD}\bigr]=\sigma_\psi\star[\psi]. \]
Suppose that $\psi'$ is another oriented optimal embedding of $\cO_K$ into $R_0^D(M)$ and write $[\psi']=\sigma\star[\psi]$ for a suitable $\sigma\in G$. Then
\begin{equation} \label{square-eq1}
\begin{split}
\bigl(\sigma_\psi\sigma^{-2}\bigr)\star[\psi']=\bigl(\sigma_\psi\sigma^{-2}\bigr)\star(\sigma\star[\psi])&=\sigma_\psi\star\bigl(\sigma^{-1}\star[\psi]\bigr)\\
&=\sigma^{-1}\star(\sigma_\psi\star[\psi])=\sigma^{-1}\star\bigl[\omega_{MD}^{-1}(\psi\circ\tau)\omega_{MD}\bigr].
\end{split}
\end{equation}
On the other hand, $\sigma_{\psi'}$ is characterized as the unique element of $G$ such that
\begin{equation} \label{square-eq2} 
\bigl[\omega_{MD}^{-1}\bigl((\sigma\star\psi)\circ\tau\bigr)\omega_{MD}\bigr]=\bigl[\omega_{MD}^{-1}(\psi'\circ\tau)\omega_{MD}\bigr]=\sigma_{\psi'}\star[\psi']. 
\end{equation}
If $\fr a\subset\cO_K$ is an ideal whose class $[\fr a]\in\Pic^+(\cO_K)$ corresponds to $\sigma$ under the reciprocity map then $\sigma^{-1}$ is represented by $\bar{\fr a}:=\tau(\fr a)$. With notation as in \S \ref{Galois-emb-subsec}, it follows that
\[ R_0^D(M)(\psi\circ\tau)(\bar{\fr a})=R_0^D(M)\psi(\fr a)=R^D_0(M)a, \]
and of course
$\bigl(a\psi a^{-1}\bigr)\circ\tau=a(\psi\circ\tau)a^{-1}$. 
Hence
\begin{equation} \label{square-eq3} 
\sigma^{-1}\star\bigl[\omega_{MD}^{-1}(\psi\circ\tau)\omega_{MD}\bigr]=\bigl[\omega_{MD}^{-1}\bigl((\sigma\star\psi)\circ\tau\bigr)\omega_{MD}\bigr], 
\end{equation}
and by comparing \eqref{square-eq1} and \eqref{square-eq2} we conclude that
\begin{equation} \label{square-eq4} 
\sigma_{\psi'}=\sigma_\psi\sigma^{-2}. 
\end{equation}
The lemma is proved. \end{proof} 
  
\begin{remark}
Equality \eqref{square-eq3} is a manifestation of the fact that the Galois extension $H^+_K/\Q$ is generalized dihedral, which is expressed by the formula $\tau\sigma\tau=\sigma^{-1}$ for all $\sigma\in G$.
\end{remark} 

Now we want to give a description of the quaternion algebra $B$ in terms of the Galois element $\sigma_\psi$ introduced in \eqref{sigma-psi-eq} (or, rather, of $\sigma_K\sigma_\psi$). This represents the analogue in our real quadratic setting of results for CM points on Shimura curves obtained by Jordan in \cite{J} (see \cite[Lemma 5.10]{GoRo} for a detailed proof of a generalization of Jordan's results).

Recall the element $\omega_\infty\in R_0^D(M)$ of norm $-1$ and consider
\[ \mathfrak D_K\star[\psi\circ\tau]=\bigl[\omega_\infty(\psi\circ\tau)\omega_\infty^{-1}\bigr]=[(\psi\circ\tau)^*]. \]
Since $\omega_\infty(\psi\circ\tau)\omega_\infty^{-1}$ has the same orientation as $\psi\circ\tau$ at all the primes dividing $MD$, it follows from Lemma \ref{galois-transitivity-lemma} that there exists $\sigma_{\fr a}\in G$ corresponding to (the class of) an ideal $\fr a$ of $\cO_K$ via the reciprocity map of class field theory such that
\begin{equation} \label{claim-eq1}
[(\psi\circ\tau)^*]=\bigl[\omega_{MD}(\sigma_{\fr a}\star\psi)\omega_{MD}^{-1}\bigr].
\end{equation}

\begin{proposition} \label{atkin-galois-prop}
There is an isomophism 
\begin{equation} \label{quaternion-algebras-isom}
B\simeq\bigg(\frac{d_K,MD{\mathrm N}_{K/\Q}(\fr a)}{\Q}\bigg)
\end{equation}
of quaternion algebras over $\Q$.
Moreover, if $\fr b$ is an ideal of $\cO_K$ satisfying \eqref{quaternion-algebras-isom} and $\sigma_{\fr b}\in G$ corresponds to the class of $\fr b$ in $\Pic^+(\cO_K)$ then the images of $\sigma_{\fr a}$ and $\sigma_{\fr b}$ in $G/G^2$ are equal.
\end{proposition}

\begin{proof} Set $i:=\psi(\sqrt{d_K})\in R_0^D(M)$. We first show the following\\

{\bf Claim.} \emph{There exists $j\in R_0^D(M)$ such that $j^2=MD{\mathrm N}_{K/\Q}(\fr a)$ and $ij=-ji$.}

\begin{proof}[Proof of Claim.] As in \S \ref{Galois-emb-subsec}, let $a\in R^D_0(M)$ be an element of positive reduced norm such that
$R_0^D(M)\psi(\fr a)=R^D_0(M)a.$ 
In particular, the reduced norm of $a$ is equal to ${\mathrm N}_{K/\Q}(\fr a)$. Then equality \eqref{claim-eq1} means that there exists $\alpha\in\Gamma_0^D(M)$ such that
\[ \alpha\omega_\infty(\psi\circ\tau)\omega_\infty^{-1}\alpha^{-1}=\omega_{MD}a\psi a^{-1}\omega_{MD}^{-1}. \]
Define $j:=\omega_\infty^{-1}\alpha^{-1}\omega_{MD}a\in R_0^D(M)$, so that
\begin{equation} \label{claim-eq2} 
\psi\circ\tau=j\psi j^{-1}. 
\end{equation}
Now we want to prove that $j^2\in\Q^\times$ and $ij=-ji$. To begin with, the relation $ij=-ji$ follows immediately from \eqref{claim-eq2} and the fact that $(\psi\circ\tau)(\sqrt{d_K})=-\psi(\sqrt{d_K})$. As for the first assertion, we use the canonical injection $B\hookrightarrow\M_2(F)$ and denote by the same symbol $\psi$ the injection $K\hookrightarrow\M_2(F)$ obtained by composing $\psi$ with this map. Recall that, in our notation, $\psi(\sqrt{d_K})=\smallmat abc{-a}\in \M_2(F)$. Write $j=\smallmat xyzw\in\M_2(F)$. Since $ij=-ji$, we see that $w=-x$, and then it follows that $j^2=\smallmat{-{\rm norm}(j)}00{-{\rm norm}(j)}$. Thus $j^2\in\Q^\times$ because ${\rm norm}(j)\in\Q$. Finally, since ${\rm norm}(j)=-MD{\mathrm N}_{K/\Q}(\fr a)\in\Z$ we conclude that $j^2=MD{\mathrm N}_{K/\Q}(\fr a)$. \end{proof}

Since $i,j\in R_0^D(M)$, in light of the above Claim we conclude that
\[ B\simeq\Q\oplus\Q i\oplus\Q j\oplus\Q ij=\bigg(\frac{d_K,MD{\mathrm N}_{K/\Q}(\fr a)}{\Q}\bigg). \]
To complete the proof of the proposition, one can proceed as in \cite[Remark 5.11]{GoRo} and note that, by genus theory, if the ideal $\fr b$ of $\cO_K$ satisfies \eqref{quaternion-algebras-isom} then $\sigma_{\fr a}$ and $\sigma_{\fr b}$ have the same image in $G/G^2$. \end{proof} 

With notation as before, it follows that
\begin{equation} \label{sigma-psi-eq2}
\sigma_\psi=\sigma_K\sigma_{\fr a}. 
\end{equation}
Now let $\lambda$ be a prime number that splits in $K$ and is such that $\lambda\equiv MD\pmod{d_K}$. A calculation with local Hilbert symbols analogous to the one in \cite[\S 5.1]{Po} shows that there is an isomorphism
\begin{equation} \label{popa-quat-eq}
B\simeq\bigg(\frac{d_K,MD\lambda}{\Q}\bigg)
\end{equation}
of quaternion algebras over $\Q$ (here we use the fact that $M$ is square-free, although the condition that $M$ be a product of odd powers of distinct primes would suffice). Choose an ideal $\fr b$ of $\cO_K$ such that ${\mathrm N}_{K/\Q}(\fr b)=\lambda$ (such a $\fr b$ exists because $\lambda$ splits in $K$). Combining isomorphism \eqref{popa-quat-eq} with Proposition \ref{atkin-galois-prop} and equality \eqref{sigma-psi-eq2}, it follows that the class of $\sigma_K\sigma_{\fr b}$ in $G/G^2$ coincides with the class $\bar\sigma$ of $\sigma_\psi$. 

Suppose that the genus character $\chi$ is associated with the pair $(\chi_1,\chi_2)$ of quadratic characters. It is now straightforward to prove the main result of this subsection.

\begin{proposition} \label{lemma4.18}
$\chi(\bar\sigma)=\chi_1(-MD)$. 
\end{proposition} 

\begin{proof} With notation as before, the congruence ${\mathrm N}_{K/\Q}(\fr b)\equiv MD\pmod{d_K}$ implies that ${\mathrm N}_{K/\Q}(\fr b)\equiv MD\pmod{d_1}$. Hence, since $\chi(\sigma_K)=\text{sign}(d_1)=\chi_1(-1)$ (cf. \cite[Remark 14.47]{Cohn}), we get that
\[ \chi(\bar\sigma)=\chi(\sigma_K)\cdot\chi(\sigma_{\fr b})=\chi_1(-MD), \]
as was to be shown. \end{proof}

\begin{remark}
Throughout this article we assume that $D>1$, i.e., that $B$ is \emph{not} isomorphic to the split algebra $\M_2(\Q)$. However, a moment's thought reveals that the arguments in this subsection remain valid also when $D=1$. With this in mind, it turns out that if $D=1$ then Proposition \ref{lemma4.18} reduces to \cite[Proposition 1.8]{BD-annals}, and so our arguments provide an alternative (and more conceptual) proof of this result.
\end{remark} 

\subsection{Derivatives of $p$-adic $L$-functions and Darmon points} \label{sec4.4}

Let $\log_E$ be as in Corollary \ref{coro}, where $\psi$ is an optimal embedding of $\mathcal O_K$ into $R_0^D(M)$. Set $J_\psi^\pm:=\log_E(P_\psi^\pm)$ and for a genus character $\chi$ of $K$ of sign $-\epsilon\in\{\pm\}$ define  
\[ J_\chi:=\sum_{\sigma\in G_K^+}\chi(\sigma)J_{\psi_\sigma}^\epsilon=\log_E(P_\chi) \] 
with $P_\chi:=\sum_{\sigma\in G_K^+}\chi(\sigma)P^\epsilon_{\psi_\sigma}\in E(K_p)$. Finally, recall the integers $n$ and $t$ in Corollary \ref{coro}. 

\begin{theorem} \label{thm-derivatives}
With notation as before, there is an equality
\[ \frac{d}{dk}\mathcal L_p(f_\infty/K,\chi,k)_{k=2}=\frac{w_{MD}\chi_1(-MD)-1}{2nt}J_\chi. \]
\end{theorem} 

\begin{proof} Recall the choice $\log_q$ of the branch of the logarithm made in \S\ref{sec2.2}. The integrals below do not depend on this choice, as observed before Definition \ref{defi-partial}. However, we will write $\log_q$ instead of a generic $\log$ to directly obtain the result to be shown. First, we have
\[ \frac{d}{dk}\mathcal L^\epsilon_p(f_\infty,\psi,k)_{k=2}=\int_\X\log_q\bigl({Q_\psi}(x,y)\bigr)d{\tilde\mu}^\epsilon_{f,\gamma_\psi}. \]
Since $Q_\psi(x,y)=c(x-z_\psi y)(x-\bar z_\psi y)$, we can write 
\[ \log_q\bigl(Q_\psi(x,y)\bigr)=\log_q(c)+\log_q(x-z_\psi y)+\log_q(x-\bar z_\psi y). \]
Since 
\[ \int_\X \log_q(c)d\tilde\mu^\epsilon_{f,\gamma_\psi}=\int_{\PP^1(\Q_p)}\log_q(c)d\pi_*(\tilde\mu^\epsilon_f)_{\gamma_\psi}=0 \]
(because $\pi_*(\tilde\mu_f)$ is cohomologous to $\mu_f$ and thus has total mass $0$ at every $\gamma$), we find that 
\[ \frac{d}{dk}\mathcal L^\epsilon_p(f_\infty,\psi,k)_{k=2}=\frac{1}{2}\int_\X\log_q(x-z_\psi y)d\mu^\epsilon_{f,\gamma_\psi}+\frac{1}{2}\int_\X\log_q(x-\bar z_\psi y)d\mu^\epsilon_{f,\gamma_\psi}. \] 
Thanks to Corollary \ref{coro}, the first term on the right hand side is equal to $-J^\epsilon_\psi/2nt$. By the same result, the second term on the right hand side is equal to $-\tau(J^\epsilon_{\psi})/2nt$. Proposition \ref{prop2.9} then shows that 
\begin{equation} \label{J-w-eq}
\frac{d}{dk}\mathcal L^\epsilon_p(f_\infty,\psi,k)_{k=2}=\frac{w_{MD}J^\epsilon_{\sigma_\psi\star\psi}-J_\psi^\epsilon}{2nt}
\end{equation}
with $\sigma_\psi$ as in \eqref{W}. Using formula \eqref{square-eq4}, it is easy to check that the map from $G_K^+$ to itself sending $\sigma$ to $\sigma_{\psi_\sigma}\sigma$ is a bijection. In light of this fact and Proposition \ref{lemma4.18}, specializing \eqref{J-w-eq} to $\psi=\psi_\sigma$ and then summing over all $\sigma\in G_K^+$ yields the result. \end{proof}

\begin{corollary} \label{coro-derivatives}
There is an equality
\[ \frac{d^2}{d^2k}L_p(f_\infty/K,\chi,k)_{k=2}=\begin{cases}J_\chi^2/{(nt)^2}&\text{if $\chi_1(-MD)=-w_{MD}$};\\[2mm]0&\text{if $\chi_1(-MD)=w_{MD}$.}\end{cases} \]
\end{corollary}

\begin{proof} The second derivative of $L_p(f_\infty/K,\chi,k)$ evaluated at $k=2$ is equal to the sum 
\[ 2\bigg(\frac{d}{dk}\mathcal L_p(f_\infty/K,\chi,k)_{k=2}\bigg)^{\!\!2}+2\mathcal L_p(f_\infty/K,\chi,2)\frac{d^2}{d^2k}\mathcal L_p(f_\infty/K,\chi,k)_{k=2}. \]
Since $\mathcal L_p(f_\infty/K,\chi,2)=0$ by \eqref{vanishing-eq}, the result follows from Theorem \ref{thm-derivatives}. \end{proof}

\subsection{Factorization formulas for $p$-adic $L$-functions}

For $j=1,2$ let $L_p(f_\infty,\chi_j,k,s)$ be the Mazur--Kitagawa $p$-adic $L$-function associated with $f_\infty$ and $\chi_j$; we refer to \cite[\S 1.4]{BD} for its definition and main properties. In particular, recall that, by \cite[Corollary 3.2]{BD-annals}, if $\chi_j(-1)=(-1)^{(k-2)/2}\tilde\epsilon$ with $\tilde\epsilon=\pm1$ (since $K$ is real, $\epsilon_K(-1)=1$, so $\chi_1(-1)=\chi_2(-1)$) we have
\begin{equation} \label{eq67}
L_p(f_\infty,\chi_j,k,k/2)=\lambda^{\tilde\epsilon}(k)\bigg(1-\frac{\chi_j(p)p^{(k-2)/2}}{a_p(k)}\bigg)^{\!\!2}\cdot\frac{\Gamma(k/2)\tau(\chi_j)}{(-2\pi i)^{(k -2)/2}}\cdot\frac{L(f_k^\sharp,\chi_j,k/2)}{\Omega_{f_k^\sharp}^{\tilde\epsilon}}
\end{equation}
where $\tau(\chi_j):=\sum_{a=1}^{d_j}\chi_j(a)e^{2\pi ia/d_j}$ is the Gauss sum attached to $\chi_j$, $\Omega_{f_k}^{\tilde\epsilon}$ is the complex period of $f_k^\sharp$ defined as in \cite[Proposition 1.1]{BD} and $\lambda^{\tilde\epsilon}(k)\in\C_p$ is defined as in \cite[Theorem 1.5]{BD}. Finally, recall that the sign of the functional equation of the $L$-function $L(f,\delta,s)$ of $f$ twisted by a quadratic Dirichlet character $\delta$ is equal to $-w_{MDp}\delta(-MDp)$ and that the $p$-adic $L$-function $L_p(f,\delta,s)$ has an exceptional zero at $s=1$ if and only if $\delta(p)=-w_p$. 

The following factorization result will play a crucial role in our subsequent arguments.

\begin{theorem} \label{thm-factor}  
There exist a neighbourhood $\U\subset\U_f$ of $2$ and a $p$-adic analytic function $\eta$ on $\U$ such that
\begin{enumerate}
\item $\eta(k)\neq0$ for all $k\in\U$ and $\eta(2)\in(\Q^\times)^2$;
\item for all $k\in\U$ there is a factorization
\[ L_p(f_\infty/K,\chi,k)=\eta(k)L_p(f_\infty,\chi_1,k,k/2)L_p(f_\infty,\chi_2,k,k/2). \] 
\end{enumerate}
\end{theorem} 

\begin{proof} By comparing Euler factors, we first notice that there is a factorization of complex $L$-functions 
\[ L(s,\pi_{f_k}\times\pi_\chi)=L(s,\pi_{f_k}\times\pi_{\chi_1})L(s,\pi_{f_k}\times\pi_{\chi_2}). \] 
On the other hand (see, e.g., \cite[p. 202]{Gel2}), for $j=1,2$ one has
\[ L\bigl(s,\pi_{f_k^\sharp}\times\pi_{\chi_j}\bigr)=\frac{\Gamma(s+(k-1)/2)}{(2\pi)^{s+(k-1)/2}}L\bigl(f_k^\sharp,\chi_j,s+(k-1)/2\bigr). \]
Let the sign of $\chi$ be $-\epsilon$. Then it follows from Theorem \ref{teo4.17} that for all even integers $k\geq4$ in $\U_f$ there is an equality 
\begin{equation} \label{eq66}
\begin{split}
L_p(f_\infty/K,\chi,k)=&\,\lambda^\epsilon_B(k)^2\bigg(1-\frac{p^{k-2}}{a_p(k)^2}\bigg)^{\!\!2}\prod_{\ell|D}\frac{\ell-1}{\ell+1}\cdot\frac{\sqrt{d_K}^{k-1}}{4\epsilon}\cdot\frac{\|\varphi_{f_k^\sharp}^{\rm JL}\|^2}{\|\varphi_{f_k^\sharp}\|^2}\cdot\frac{\Gamma(k/2)^2}{(-i\pi)^k}\\
&\times L\bigl(f_k^\sharp,\chi_1,k/2\bigr)L\bigl(f_k^\sharp,\chi_2,k/2\bigr).
\end{split}
\end{equation}
Since $p$ is inert in $K$, we have $\chi_1(p)=-\chi_2(p)$, and then formulas \eqref{eq67} and \eqref{eq66} imply that  
\begin{equation} \label{eq-factor+}
L_p(f_\infty/K,\chi,k)=\eta(k)L_p(f_\infty,\chi_1,k,k/2)L_p(f_\infty,\chi_2,k,k/2)
\end{equation}
for all even integers $k\geq4$, where $\eta(k)\in\C_p$ can be made explicit by comparing \eqref{eq67} and \eqref{eq66}.

Now we show that $\eta$ can be extended to an analytic function with the prescribed properties. For this we need an argument similar to the one in the proof of \cite[Proposition 5.2]{BD}. Consider the set of primitive Dirichlet characters $\chi_1'$ of conductor prime to $MDp$ such that 
\begin{itemize}
\item $\chi_1'(-1)=\chi_1(-1)=\chi_2(-1)$; 
\item $\chi_1'(-MD)=-w_{MD}$; 
\item $\chi_1'(p)=-a_p=w_p$. 
\end{itemize}
The condition $\chi_1'(-MD)=-w_{MD}$ implies that $L(f,\chi_1',s)$ vanishes of even order at $s=1$, while the condition $\chi_1'(p)=-a_p=w_p$ implies that $L_p(f,\chi_1',1)$ does not have an exceptional zero, so that $L_p(f,\chi_1',1)\in\Q^\times$. Fix a character $\chi'_1$ as above such that $L(f,\chi'_1,1)\not=0$ (by the main result of \cite{MM}, there are infinitely many).

Consider now the set of primitive Dirichlet characters $\chi_2'$ of conductor prime to $MDp$ and the conductor of $\chi_1'$ satisfying
\begin{itemize} 
\item $\chi_2'(-1)=\chi_1'(-1)$;
\item $\chi_2'(\ell)=\chi_1'(\ell)$ for $\ell|pM$;
\item $\chi_2'(\ell)=-\chi_1'(\ell)$ for $\ell|D$.
\end{itemize} 
Since $D$ has an even number of prime factors, $\chi_2'(-MD)=\chi_1'(-MD)$ and $L(f,\chi_2',s)$ also vanishes to even order at $s=1$. Also, since $\chi_2'(p)=-a_p$, the $p$-adic $L$-function $L_p(f,\chi_2',s)$ does not have an exceptional zero. Hence, again by \cite{MM}, there are infinitely many characters $\chi_2'$ as above such that $L(f,\chi_2',1)\neq0$ and $L_p(f,\chi_2',1)\neq 0$. Choose a pair $(\chi_1',\chi_2')$ as prescribed above and denote by $d_1'$ and $d_2'$ the conductors of $\chi_1'$ and $\chi_2'$, respectively. Let $d_{K'}:=d_1'\cdot d_2'$ and consider the real quadratic field $K':=\Q(\sqrt{d_{K'}})$ and the genus character $\chi'$ of $K'$ corresponding to the pair $(\chi_1',\chi_2')$. In particular, by definition, the sign of $\chi'$ is equal to the sign of $\chi$. Recall that $G_{K'}^+:=\Gal(H_{K'}^+/K')$ and fix a set $\{\psi'_\sigma\}$ of representatives of the $\Gamma_0^D(Mp)$-equivalence classes of oriented optimal embeddings of $\mathcal O_{K'}$ into $R_0^D(Mp)$ as described in Remark \ref{remark4.16}, so that this family is also a set of representatives for the $\Gamma_0^D(M)$-equivalence classes of oriented optimal embeddings of $\mathcal O_{K'}$ into $R_0^D(M)$. Now consider the expression 
\begin{equation} \label{p-adic-L-split}
L_p(f_\infty/K',\chi',k):=\Bigg(\sum_{\sigma\in G_{K'}^+}\chi'(\sigma)\int_{\Z_p^\times\times\Z_p^\times}Q^{(k)}_{\psi'_\sigma}(x,y)d\tilde\mu_{f,\gamma_{\psi'_\sigma}}^{\epsilon}\Bigg)^{\!\!2}.
\end{equation}
Combining the arguments in the proof of Proposition \ref{prop4.15} with Proposition \ref{prop-split} and Corollary \ref{coro-split}, we see that 
\[ \int_{\Z_p^\times\times\Z_p^\times}d\tilde\mu_{f,\gamma_{\psi'_\sigma}}^\epsilon=2\lambda_B^\epsilon(2)I^\epsilon\bigl(f^{{\rm JL},\sharp},\tau,\gamma_{\psi'_\sigma}\bigr)(1) \]
and  
\[ \int_{\Z_p^\times\times\Z_p^\times}Q^{(k)}_{\psi'_\sigma}(x,y)d\tilde\mu_{f,\gamma_{\psi'_\sigma}}^\epsilon=\lambda_B^\epsilon(k)\bigg(1-\frac{w_pp^{(k-2)/2}}{a_p(k)}\bigg)^{\!\!2}I^\epsilon\bigl(f^{{\rm JL},\sharp}_k,\tau,\gamma_{\psi'_\sigma}\bigr)\Big(Q_{\psi'_\sigma}^{(k)}(x,y)\Big) \]
for $f_k$ with trivial character. Set $f_2^\sharp:=f$, $f_2^{\rm JL,\sharp}:=f^{\rm JL}$ and $\varphi_{f_2}^{\rm JL,\sharp}:=\varphi_f^{\rm JL}/2$. For all even $k\geq2$ such that $f_k$ has trivial character we obtain the interpolation formula  
\[ \begin{split}
   L_p(f_\infty/K',\chi',k)=&\,\lambda_B^\epsilon(k)^2\bigg(1-\frac{w_pp^{(k -2)/2}}{a_p(k)}\bigg)^{\!\!k}\prod_{\ell|D}\frac{\ell-1}{\ell+1}\cdot\frac{(2i)^k\sqrt{d_K}^{k-1}}{4\epsilon}\\
   &\times\frac{\|\varphi_{f_k}^{\rm JL,\sharp}\|^2}{\|\varphi_{f_k^\sharp}\|^2}\cdot L_K(1/2,\pi_{f_k^\sharp}\times\pi_\chi).
   \end{split} \]
Therefore we get the factorization 
\[ L_p(f_\infty/K',\chi',k)=\eta(k)L_p(f_\infty,\chi'_1,k,k/2)L_p(f_\infty,\chi'_2,k,k/2) \]
where $\eta$ is the function appearing in \eqref{eq-factor+}.
 
With these choices in force, $L_p(f_\infty/K',\chi',k)$, $L_p(f_\infty,\chi'_1,k,k/2)$ and $L_p(f_\infty,\chi'_2,k,k/2)$ do not vanish at $k=2$, so there exists a neighbourhood $\U$ of $2$ such that they do not vanish at any point of $\U$. Recall that $\lambda^{\tilde\epsilon}(2)$ is chosen to be equal to $1$. Thus the quotient 
\[ \eta(k,\chi_1',\chi_2'):= \frac{L_p(f_\infty/K',\chi',k)}{L_p(f_\infty,\chi'_1,k,k/2)L_p(f_\infty,\chi'_2,k,k/2)} \]
is an analytic function on $\U$ which does not vanish at any point of $\U$, and we define $\eta(k):=\eta(k,\chi_1',\chi_2')$. (Observe that $\eta(k,\chi'_1,\chi'_2)$ does not depend on the choice of $\chi'_1$ and $\chi'_2$ made above: any other choice defines a function which coincides with $\eta(k,\chi'_1,\chi'_2)$ on the dense subset of points $k\in\U$ with trivial character because, for these points, $\eta(k)$ does not depend on $\chi'_1$, $\chi'_2$.)

To conclude the proof, we need only show the statement about $\eta(2)$. For this, we use the following argument. Fix a prime $d|D$ and consider the set of Dirichlet characters $\chi_3'$ and $\chi_4'$ of conductor prime to $MDp$ satisfying the following conditions: 
\begin{itemize}
\item $\chi_1'(-1)=\chi_2'(-1)=-\chi_3'(-1)=\chi_4'(-1)$;
\item $\chi_1'(\ell)=\chi_2'(\ell)=\chi_3'(\ell)=\chi_4'(\ell)$ for all primes $\ell|pM$; 
\item $\chi_1'(d)=-\chi_2'(d)=-\chi_3'(d)=-\chi_4'(d)$;
\item $\chi_1'(\ell)=-\chi_2'(\ell)=\chi_3'(\ell)=-\chi_4'(\ell)$ for all primes $\ell|(D/d)$.
\end{itemize}
To simplify notations, set $\Lambda_j:=L_p(f_\infty,\chi'_j,2,1)$ for $j=1,\dots,4$. Note that $\chi_3'(-MD)=\chi_4'(-MD)=-w_{MD}$ and $\chi_3'(p)=\chi_4'(p)=w_p$. Thus, as above, the order of vanishing of $L(f,\chi_3',s)$ and $L(f,\chi'_4,s)$ at $s=1$ is even and $L_p(f,\chi_3',s)$ and $L_p(f,\chi_4',s)$ do not have an exceptional zero at $s=1$. Using \cite{MM} again, we choose $\chi'_3$ and $\chi'_4$ so that $\Lambda_3\not=0$ and $\Lambda_4\neq 0$. 

For the pairs of indices $(1,3)$, $(3,4)$ and $(2,4)$ we consider quadratic fields $K'_{i,j}$ and genus characters $\chi'_{i,j}$ of $K'_{i,j}$ such that the associated Dirichlet characters are $(\chi'_i,\chi'_j)$; in particular, if $d_j'$ is the conductor of $\chi_j'$ then the discriminant of $K_{i,j}'$ is $d_i'd_j'$. Note that
\begin{itemize}
\item $K'_{1,3}$ and $K'_{3,4}$ are imaginary, while $K'_{2,4}$ is real; 
\item the primes dividing $Mp$ are split in $K'_{1,3}$ and $K'_{3,4}$ and the number of primes dividing $D$ which are inert in these fields is odd;
\item the primes dividing $MDp$ are split in $K'_{2,4}$. 
\end{itemize} 
Thanks to \cite[Proposition 5.1]{BD}, we have 
\[ L_p(f_\infty/K'_{i,j},\chi'_{i,j},2)=\langle \phi_2,\phi_2\rangle\Lambda_i\Lambda_j\qquad\text{for $(i,j)=(1,3)$ and $(i,j)=(3,4)$}, \]
where we use the notation in \cite{BD}. In particular, recall that the $p$-adic $L$-function in the equation above is defined in \cite[\S 3.2]{BD} and that $\langle\phi_2,\phi_2\rangle\in\Q^\times$ is defined in \cite[p. 395]{BD}. Since, by construction, $L_p(f_\infty/K'_{i,j},\chi'_{i,j},2)$ is a square in $\Q^\times$, it follows that $\Lambda_i\Lambda_j\in\langle \phi_2,\phi_2\rangle\cdot(\Q^\times)^2$ for  $(i,j)=(1,3)$ and $(i,j)=(3,4)$. 

By \cite[Theorem 3.6]{BD-annals}, there is also an equality 
\[ L_p(f_\infty/K'_{2,4},\chi'_{2,4},2)=\Lambda_2\Lambda_4 \]
where, this time, $L_p(f_\infty/K'_{2,4},\chi'_{2,4},2)$ is defined in \cite[Definition 3.4]{BD-annals}. Again, by construction, $L_p(f_\infty/K'_{2,4},\chi'_{2,4},2)$ is a square in $\Q^\times$, so the same is true of the product $\Lambda _2\Lambda_4$. Now the factorization
\[ \prod_{i=1}^4\Lambda_i=(\Lambda_1\Lambda_2)\cdot(\Lambda_3\Lambda_4)=(\Lambda_1\Lambda_3)\cdot(\Lambda_2\Lambda_4)\in(\Q^\times)^2 \] 
shows that $\Lambda_1\Lambda_2$ is a square in $\Q^\times$. Thanks to the interpolation formula  
\[ L_p(f_\infty/K',\chi',2)=\eta(2)\Lambda_1\Lambda_2 \]
proved before and the fact that, by construction, $L_p(f_\infty/K',\chi',2)$ is a square in $\Q^\times$ (keep in mind the normalization chosen in part (1) of Proposition \ref{prop-pi-2}), it follows that the same is true of $\eta(2)$, as was to be shown. \end{proof}

\begin{remark}
The function $L_p(f_\infty/K',\chi',k)$ appearing in \eqref{p-adic-L-split} can, of course, be interpreted as a $p$-adic $L$-function attached to the Hida family $f_\infty$, the real quadratic field $K'$ and the genus character $\chi'$. In particular, observe the following facts.

1) Since $p$ is split in $K'$, the set $\Z_p^\times\times\Z_p^\times$ is stable under the action of the elements $\gamma_{\psi_\sigma}$, hence, thanks to Lemma \ref{well-def} and the density argument in the proof of Proposition \ref{prop2.24}, the function $L_p(f_\infty/K',\chi',k)$ is well defined independently of the choices of the representative $\tilde\mu_f^\epsilon$ and of the set $\{\psi_\sigma\}_{\sigma\in G_{K'}^+}$. 

2) As in Remark \ref{rem4.24}, one can show that for all $k\in\U_f$ there is an equality 
\[ \int_{\Z_p^\times\times\Z_p^\times} Q_\psi^{(k)}(x,y)d\tilde\mu^\pm_{f,\gamma_\psi}=\int_{\Z_p^\times\times\Z_p^\times}\langle Q_\psi(x,y)\rangle^{\frac{k-2}{2}}d\tilde\mu^\pm_{f,\gamma_\psi}. \] 
To show this, one proceeds as in the proof of \cite[Lemma 3.7]{BD} upon noticing that $\psi$, $\Z_p\times\Z_p$ and $\Z_p^\times\times\Z_p^\times$ correspond to $\Psi$, $L_\Psi$ and $L_\Psi''$ in the notations of \cite{BD}. Since we need the function in \eqref{p-adic-L-split} only as an auxiliary object, we will not develop the theory of $p$-adic $L$-functions in this context any further. 
\end{remark}

\section{Rationality results for Darmon points on elliptic curves} 

Set $N:=MDp$, let $E_{/\Q}$ be an elliptic curve of conductor $N$ and let $K$ be a real quadratic field as in \S \ref{sec2.1}. Recall that $D>1$ is a square-free product of an \emph{even} number of primes.  

\begin{theorem} \label{theorem5.1}
Let $\chi$ be a genus character of $K$ corresponding to a pair $(\chi_1,\chi_2)$ such that $\chi_i(-MD)=-w_{MD}$ for $i=1,2$. 
\begin{enumerate}
\item There exists a point ${\bf P}_\chi\in E(H_\chi)^\chi$ and a rational number $c\in\Q^\times$ such that 
\[ J_\chi=c\log_E({\bf P}_\chi). \]
\item The point ${\bf P}_\chi$ is of infinite order if and only if $L^\prime(E/K,\chi,1)\neq 0$.
\item A suitable integral multiple of $P_\chi$ belongs to the natural image of $E(H_\chi)^\chi$ in $E(K_p)$. In particular, $P_\chi$ coincides with the image of a global point in $E(K_p)\otimes_\Z\Q$.  
\end{enumerate} 
\end{theorem} 

\begin{proof} As in the proof of \cite[Theorem 4.3]{BD-annals}, order $\chi_1,\chi_2$ in such a way that ${\rm sign}(E,\chi_1)=-1$. Then $\chi_1(p)=-w_p=a_p$ and, thanks to \cite[Theorem 5.4]{BD}, $L_p(f_\infty,\chi_1,k,k/2)$ vanishes of order at least $2$ at $k=2$. Furthermore, there are a global point ${\bf P}_{\chi_1}\in E\bigl(\Q(\sqrt{d_1})\bigr)^{\chi_1}$, which is of infinite order if and only if $L^\prime(E/\Q,\chi_1,1)\neq0$, and an integer $s$ such that
\begin{equation} \label{eq-der}
\frac{d^2}{dk^2}L_p(f_\infty,\chi_1,k,k/2)_{k=2}=s\log_E({\bf P}_{\chi_1})^2
\end{equation}
and
\begin{equation} \label{eq-n}
s\equiv L^*(f,\psi,1)\pmod{(\Q^\times)^2}
\end{equation}
for \emph{any} primitive Dirichlet character $\psi$ for which $L(f,\psi,1)\neq0$, $\psi(p)=-\chi(p)$ and $\psi(\ell)=\chi(\ell)$ for all primes $\ell|MD$. Here the \emph{algebraic part} $L^*(f,\psi,1)$ of $L(f,\psi,1)$ is as in \cite[eq. (24)]{BD-annals}.

From Theorem \ref{thm-factor} and the fact that $L_p(f_\infty,\chi_1,k,k/2)$ has order of vanishing at $k=2$ greater than or equal to $2$, it follows that 
\begin{equation} \label{eq-der-fac}
\frac{d^2}{dk^2}L_p(f_\infty/K,\chi,k)_{k=2}=\eta(2)\frac{d^2}{dk^2}L_p(f_\infty,\chi_1,k,k/2)_{k=2}L_p(f_\infty,\chi_2,2,1).
\end{equation}
First suppose that $L'(E/K,\chi,1)\neq0$. In this case $L_p(f_\infty,\chi_2,2,1)=2L^*(f,\chi_2,1)\in\Q^\times$. Set $m:=L^*(f,\chi_2,1)$. By \eqref{eq-n}, there exists $r\in\Q^\times$ such that $r^2=s/m$. Define ${\bf P}_\chi:={\bf P}_{\chi_1}$ in this situation. Setting $v:=ntrm\in\Q^\times$, a combination of \eqref{eq-der}, \eqref{eq-der-fac} and Corollary \ref{coro-derivatives} yields 
\[ J_\chi^2=\eta(2)(nt)^2sm\log_E({\bf P}_{\chi_1})^2=\bigl(v\sqrt{\eta(2)}\bigr)^2\log_E({\bf P}_\chi)^2, \] 
from which we obtain the desired result for $c:=\pm v\sqrt{\eta(2)}\in\Q^\times$ (recall that, by Theorem \ref{thm-factor}, $\eta({2})$ is a square in $\Q^\times$). 

Now suppose that $L'(E/K,\chi,1)=0$. In this case $L_p(f_\infty,\chi_2,2,1)=2L^*(f,\chi_2,1)=0$, and the result is an immediate consequence of \eqref{eq-der-fac} and Corollary \ref{coro-derivatives} on setting ${\bf P}_\chi:=0$. This completes the proof of parts (1) and (2). 

Finally, since $J_\chi=\log_E(P_\chi)$, if $c=a/b$ with $a,b\in\Z$ then part (2) gives the equality 
\[ \log_E(bP_\chi)=\log_E(a{\bf P}_\chi), \]
which implies part (3) since the kernel of $\log_E$ is the torsion subgroup of $E(K_p)$. \end{proof}

\end{document}